 \documentclass{amsart}
\usepackage{graphicx}

 \usepackage{amscd}
                     \usepackage{amssymb}
 \newcommand{\be}{\begin{equation}}
       \newcommand{\ee}{\end{equation}}
       \newcommand{\ba}{\begin{eqnarray}}
        \newcommand{\ea}{\end{eqnarray}}
 \newcommand{\ban}{\begin{eqnarray*}}
 \newcommand{\ean}{\end{eqnarray*}}

 \newcommand{\ra}{\rightarrow}

  \newcommand{\Pf}{\noindent {\bf Proof:} }

 \newcommand{\sect}[1]{\section{#1} \setcounter{equation}{0}}

 \newtheorem{theo}{Theorem}[section]

\begin{document}
 \newtheorem{example}[theo]{Example}
 \newtheorem{defn}[theo]{Definition}
 \newtheorem{ques}[theo]{Question}
 \newtheorem{lem}[theo]{Lemma}
 \newtheorem{lemma}[theo]{Lemma}
 \newtheorem{prop}[theo]{Proposition}
 \newtheorem{coro}[theo]{Corollary}
 \newtheorem{ex}[theo]{Example}
 \newtheorem{note}[theo]{Note}
 \newtheorem{conj}[theo]{Conjecture}
 \newtheorem{rmk}[theo]{Remark}

 \newcommand{\inj}{\mbox{inj}}
 \newcommand{\vol}{\mbox{vol}}
 \newcommand{\diam}{\mbox{diam}}
 \newcommand{\Ric}{\mbox{Ric}}
 \newcommand{\Iso}{\mbox{Iso}}
 \newcommand{\Hess}{\mbox{Hess}}
 \newcommand{\divg}{\mbox{div}}

\title{The Cut-off Covering Spectrum}
\author{Christina Sormani}
\email{sormanic@member.ams.org}
 \thanks{Partially supported by
a grant from the City University of New York PSC-CUNY Research Award
Program.}
\address{CUNY Graduate Center and Lehman College}
\author{Guofang Wei}
 \email{wei@math.ucsb.edu}
 \address{Department of Mathematics, UCSB, Santa Barbara, CA 93106}
\thanks{Partially
 supported by NSF Grant \# DMS-0505733.}
 \date{}

 \begin{abstract}
We introduce the $R$ cut-off covering spectrum and the cut-off
covering spectrum of a metric space or Riemannian manifold.
The spectra measure the sizes of localized holes in the space and
are defined using covering spaces called $\delta$ covers and $R$
cut-off $\delta$ covers. They are investigated using $\delta$
homotopies which are homotopies via grids whose squares are mapped
into balls of radius $\delta$.

On locally compact spaces, we prove that these new spectra are
subsets of the closure of the length spectrum.   We prove the $R$
cut-off covering spectrum is almost continuous with respect to the
pointed Gromov-Hausdorff convergence of spaces and that the cut-off
covering spectrum is also relatively well behaved. This is not true
of the covering spectrum defined in our earlier work which was shown
to be well behaved on compact spaces. We close by analyzing these
spectra on Riemannian manifolds with lower bounds on their sectional
and Ricci curvature and their limit spaces.
  \end{abstract}

 \maketitle

\sect{Introduction} \label{sect1}

Metric spaces and Riemannian manifolds
are often studied using Gromov Hausdorff convergence and Gromov's
compactness theorem. However, this convergence, reviewed in
Section~\ref{Sect5}, does not preserve the topology of the space.  Thinner and
thinner flat tori converge to circles, thus losing a generator of
the fundamental group. Sequences of surfaces of higher and higher
genus can converge to the Hawaii Ring, a space with an infinitely
generated fundamental group and no universal cover
[Example~\ref{Hawaii}].  Sequences of capped cylinders can be
seen to converge in the pointed Gromov-Hausdorff sense to
cylinders if the base points slide out to infinity [Example~\ref{cylinder}].

Adding curvature conditions to the spaces in question both restrict
their topology and topology of the limit spaces.  Cheeger-Gromoll's Soul
Theorem not only demonstrates that complete noncompact spaces with
nonnegative sectional curvature have finite topological type, but also that
their ``holes" are located in a compact soul \cite{Cheeger-Gromoll}.
Perelman proved a geometric extension of this result for their limit spaces
using work of Sharafutdinov \cite{Sharafutdinov}\cite{Perelman1}.  More
recently Cheeger-Colding
have proven a number of results concerning the limits of manifolds with
nonnegative Ricci curvature (c.f. \cite{Cheeger} and \cite{Wei-survey}).
The topology of such spaces has been studied extensively by a number of
mathematicians (c.f. \cite{ShSo2}).  As one examines this work, it becomes
clear
that it is not only of importance to understand the topological question
concerning the existence of holes in these spaces but also to examine the
geometric properties of these holes.

In \cite{SoWei3}, the authors defined the covering spectrum of a
compact length space, $K$.  This spectrum measures
the size of the one dimensional holes in the space and is closely
related to the length spectrum: every element
in the covering spectrum is half the length of a closed geodesic,
\be \label{oldlengththm}
CovSpec(K)\subset (1/2)Length(K).
\ee
The covering spectrum is
empty when the space is simply connected or is its own universal
cover.  It is determined using a sequence of covering spaces called
$\delta$ {\em covers} which unravel curves that do not fit in balls
of radius $\delta$.   We proved that when compact length spaces
$K_i$ converge in the Gromov-Hausdorff sense to a compact length
space $K$, then their covering spectra converge in the Hausdorff
sense:
\be
d_{H}(CovSpec(K_i)\cup \{0\}, CovSpec(K)\cup \{0\})\to 0.
\ee
It is possible for elements to converge to $0$ as they do on the
sequence of thinner and thinner tori, at which point they disappear
and are no longer in the covering spectrum. However, elements which
converge to a positive value do not disappear in the limit.
Furthermore, an element in $CovSpec(K)$
is the limit of elements in $CovSpec(K_i)$.
In particular, the covering spectrum of the limit space, $K$, of a
sequence of simply connected spaces, $K_i$, is empty \cite{SoWei3}.

When studying complete locally compact spaces, it is natural to employ
{\em pointed Gromov-Hausdorff convergence}.  The covering spectrum is not
continuous with respect to this convergence.  Sequences of
manifolds, $X_i$ with handles sliding out to infinity converge to a
space $X$ with no handles, so that we can have $\delta\in
CovSpec(X_i)\, \forall i\in \mathbb{N}$ yet $CovSpec(X)=\emptyset$
[Example~\ref{disappear}].  It
is even possible for there to be an element in the covering spectrum
of the limit space when $CovSpec(X_i)=\emptyset\,\forall i \in
\Bbb{N}$ [Example~\ref{ex5.2}].  These difficulties arise because the
pointed Gromov-Hausdorff convergence is defined as the
Gromov-Hausdorff limit of balls of radius $R$ where the convergence
can be slower as we take larger values of $R$.

Further difficulties are caused by the lack of compactness on a
single space. Even on a locally compact space
the covering spectrum is no longer closely related to
the length spectrum on a noncompact space: there can be holes which
extend to infinity and decrease in size [Examples~\ref{notslip}
and~\ref{slipex}].
Those that decrease to $0$ are not detected by the covering spectrum
and those that decrease to a constant cause an element in the
covering spectrum to exist which is not $1/2$ the length of a closed
geodesic.  This will be explored in a future paper.  Here we define a
new spectrum which resolves many of these difficulties.

In this paper we introduce the $R$ {\em cut-off covering spectrum}
and the {\em cut-off covering spectrum} to overcome these
difficulties. The $R$ cut-off covering spectrum of a pointed space
$(X,x)$ detects holes which do not extend outside the closed ball
$\bar{B}_x(R)$.  The cut-off covering spectrum detects holes
which do not extend to infinity. A cylinder only has a hole
extending to infinity, so its cut-off covering spectrum is empty.
We prove that on a complete locally compact length space, $X$,
 both of these spectra are contained in the closure of the length
spectrum because the holes they detect are localized
[Theorem~\ref{properRcutinL}
and Corollary~\ref{coro4.16}].  Local compactness is seen to be necessary
in Example~\ref{boundednotinL}.

We prove that the $R$ cut-off covering spectrum is continuous with
respect to the pointed Gromov-Hausdorff convergence of the locally
compact spaces [Theorem~\ref{thm5.2}].  This result is not an immediate
extension of our compact results because the $R$ cut-off spectrum is
not uniformly localized: it detects any hole which passes into
$\bar{B}_x(R)$ no matter how far out part of the hole extends.  While the
elements of the covering spectrum of a compact space are bounded
above by the diameter of the space, there is no upper bound on an
element in the $R$ cut-off covering spectrum [Example~\ref{ex5.5}].
In Example~\ref{ex5.5}, one sees a sequence of $X_i$ with increasingly
large holes such that the hole snaps open to give a simply
connected limit $X$.  One aspect of our
theorem says that if a sequence of spaces $(X_i,x_i)$ have elements
\be \delta_i\in CovSpec_{cut}^R(X_i,x_i) \ee which diverge to
infinity, then the holes they detect always snap open in the limit and are
no longer holes at all.

Another difficulty in our noncompact setting arises from
the fact that the $R$ cut-off
covering spectrum is defined using covering spaces and, as such,
homotopies which extend far outside  $\bar{B}_x(R)$ could influence
the value of $CovSpec_{cut}^R(X,x)$.
To handle this issue we develop the concept of the $\delta$ {\em
homotopy} first introduced \cite{SoWei1}.  A closed curve is $\delta$
homotopic to a point if it lifts as a closed curve to the $\delta$
cover of the space.  We introduce $\delta$ homotopies: maps from
rectangular grids to the space which map squares into balls of
radius $\delta$ [Lemma~\ref{lemma3.3}]. This allows us to control
the location of
the maps and, in particular, we prove that if a curve is $\delta$
homotopic to a point, then it is $\delta$ homotopic in a bounded
region to a collection of possibly trivial loops lying near the
boundary of that region [Lemma~\ref{lemma3.5}].
Later we apply this to localize
subsets of the cut-off covering spectrum [Proposition~\ref{isomballs}].
We also bound the lengths of curves in a region $A$ with certain
$\delta$ homotopic properties in terms of the number of disjoint
balls of radius $\delta/5$ that fit within $A$ [Lemma~\ref{lem3.7}].  This
is useful later for uniformly bounding the size of holes which are
detected by the $R$ cut-off covering spectrum in a Gromov-Hausdorff
converging sequence of balls.

In Section~\ref{review} we review crucial concepts from prior papers,
simplifying some and clarifying others.  We begin with the classical
notions of
covering spaces and length spaces and review which metric properties 
lift from the base space to the cover in Subsection~\ref{sect2.0}.
Subsections~\ref{sect2.1}
and~\ref{sect2.2} review the definitions of delta covers and the
covering spectrum as in \cite{SoWei1} and \cite{SoWei3} respectively.
Subsection~\ref{sect2.3} provides a simplified but equivalent definition of
the covering spectrum for a complete length space with a universal
cover [Defn~\ref{groupcovspec} and Thm~\ref{balltoloop}].  
Note that we do not require the existence of a simply connected cover,
but rather just a covering space which covers all other covering spaces
in the sense of Spanier \cite{Sp}.

In Section~\ref{delta-homotopy} we develop the 
theory of $\delta$ homotopy.  Similar but distinct concepts
appear in work of Beretovskii-Plaut \cite{BP2} 
which focuses on the construction of a uniform universal cover which
is not a covering space but is simply connected.  
While the lemmas in
Section~\ref{delta-homotopy} are very intuitive and have explanatory
diagrams, the proofs are necessarily technical and may be skipped by
the reader.  

In Section~\ref{cut-off-spectrum} we introduce the cut-off
covering spectra of pointed metric spaces $(X,x)$.
We begin by defining the $R$ cut-off $\delta$ covers,
$\tilde{X}^{\delta,R}_{cut}$, which unravel curves that
are not $\delta$ homotopic to loops outside $\bar{B}_x(R)$.
We prove that when $X$ is a complete locally compact length space,
$\tilde{X}^{\delta,R}_{cut}$ have
unique limits as $R$ diverges to infinity and call these limits the
cut-off $\delta$ covers [Prop~\ref{prop4.6}].  These covers unravel
holes which do not extend to infinity. The $R$ cut-off covering spectrum,
$CovSpec^R_{cut}(X,x)$
 is
defined using the $R$ cut-off $\delta$ covers for any metric
space while the cut-off
covering spectrum $CovSpec_{cut}(X)$
is defined using the cut-off $\delta$ covers and is basepoint invariant
[Definitions~\ref{defn4.4} and~\ref{defn4.8}].  By Prop~\ref{prop4.6},
any complete locally compact length space has a well defined cut-off
covering spectrum but we believe algebraic techniques might
be used to prove it is well defined for a much larger class of
metric spaces.  

In Section~\ref{sect4.3} we relate these spectra to the covering spectra
and to each other, showing in particular that for any $R_1<R_2$
and any basepoint $x$ in a metric space $X$ with a well defined
cut off covering spectrum, we have
\be
CovSpec_{cut}^{R_1}(X,x)\subset CovSpec_{cut}^{R_2}(X,x) \subset
CovSpec_{cut}(X) \subset CovSpec(X).
\ee

In Section~\ref{sect4.4},  we prove Theorem~\ref{properRcutinL} that for
a complete
locally compact length space, $X$,  if $\delta \in CovSpec^R_{cut}(X)$
then $2\delta \in Length(X)$ which we write as
$CovSpec^R_{cut}(X)\subset (1/2)Length(X)$.
As a corollary we then show
\be \label{firstlower}
CovSpec_{cut}(X) \subset Cl_{lower}((1/2)Length(X))
\ee
where
$Cl_{lower}(A)$ is the lower semiclosure of the set $A\subset
\Bbb{R}$. The lower semiclosure is defined and explored in the
appendix, where we prove any spectrum defined in a manner similar to
these spectra are lower semiclosed sets [Theorem~\ref{appendixthm}].
Example~\ref{ex4.2} demonstrates the necessity of the lower semiclosure in
(\ref{firstlower}) .

In Section~\ref{sect4.5} we study various topological
conditions on metric spaces with well defined cut off covering spectra.
We first recall the loops to infinity property
defined in \cite{So-loops} and relate this concept
to the emptiness of the cut-off covering spectrum
[Theorem~\ref{thmloopstoinfty} and Theorem~\ref{thm-one-end}].
Corresponding examples are presented as well.
Then we describe the cut-off covering spectrum on
product spaces [Theorem~\ref{topcross}].

In Section~\ref{sectribbon}, we introduce a new construction
of length spaces which are not locally compact.  This construction
consists of attaching a ``pulled ribbon'' to a given space along a
line.  Example~\ref{ribbon1} demonstrates the necessity of
the local compactness condition in Theorem~\ref{thmoneend}.
Example~\ref{ribbon2} demonstrates that a space with
an empty length spectrum can have a nontrivial cut-off covering
spectrum demonstrating the necessity of local compactness
in  Theorem~\ref{properRcutinL}.

In Section~\ref{sect4.6} we
localize the $R$ cut-off covering spectrum using the $\delta$
homotopies as mentioned above. Proposition~\ref{isomballs} shows
subsets of the $R$ cut-off covering spectra agree on spaces
with isometric balls of sufficient size.

In Section~\ref{sect4.7} we explore
\be
CovSpec_{cut}^{R_2}(X)\setminus CovSpec_{cut}^{R_1}(X) \ \ \mbox{when} \ R_2
> R_1.
\ee
In particular Propositions~\ref{prop4.26} and~\ref{prop4.27} together imply
that these two spectra are equivalent for $R_2$ sufficiently
close to $R_1$ on locally compact spaces .

In Section~\ref{sect5} we introduce Gromov-Hausdorff convergence, first
reviewing the definitions.  In  Section~\ref{sect5.1} we provide examples
demonstrating why the covering spectrum is not continuous with
respect to pointed Gromov-Hausdorff convergence: elements can shrink
to $0$, disappear in the limit, suddenly appear in the limit, or
diverge to infinity.

In Section~\ref{sect5.2} we prove the continuity of the $R$ cut-off covering
spectrum [Theorem~\ref{thm5.2}]
and provide examples clarifying why it is necessary to
slightly change $R$ to obtain this continuity.  The
proof requires two propositions: one controlling the fundamental
groups of the $R$ cut-off $\delta$ covers and the other proving the
$R$ cut-off $\delta$ covers converge. It also strongly relies on the
results on $\delta$ homotopies and localization proven in the
earlier sections.

In Section~\ref{sect5.3} we prove Theorem~\ref{noappear} which states that
\be
\textrm{ for any } \delta \in CovSpec_{cut}(X), \textrm{ there is }
\delta_i \in CovSpec_{cut}(X_i)
\ee
  such that $\delta_i \rightarrow \delta$.
In particular if $X_i$ are simply connected locally compact spaces
that converge to a locally compact space $X$ in the pointed
Gromov-Hausdorff sense then $CovSpec_{cut}(X)=\emptyset$
[Corollary~\ref{limsimp}]. This limit space need not be simply
connected as can be seen in Example~\ref{ex5.2}.

Further directions of study are suggested in Question~\ref{questproperconv}
and Remark~\ref{rmkproperconv}.

In Section~\ref{sect5.4} we prove the pointed Gromov-Hausdorff limits of
simply connected spaces either have the loops to infinity property
or two ends [Theorem~\ref{thm5.7}].
In Section~\ref{sect5.5} we investigate the cut-off
covering spectra of tangent cones at infinity, proving
in Theorem~\ref{tannoappear} that spaces with bounded cutoff covering
spectra have tangent cones at infinity with empty covering spectra.

We close the paper with Section~\ref{sect6} on applications to spaces with
curvature bounds.  Section~\ref{sect6.1} discusses manifolds with
nonnegative sectional curvature and consequences of the Cheeger-Gromoll
Soul Theorem and work of Sharafutdinov and Perelman. Theorem~\ref{sectthm}
applies to length spaces with curvature bounded below as
well and states that if $S$ is the soul of the manifold, then
\be \label{souleqn*}
CovSpec(S)=CovSpec(T_{R}(S))=CovSpec(M)
\ee
where $T_R(S)$ is the tubular neighborhood of the soul.

In Section~\ref{sect6.2} we apply our convergence results to
obtain an almost soul theorem [Theorem~\ref{almostsoul}]
which says that locally (\ref{souleqn*}) is approximately true.
Corollary~\ref{sectge-epsi} descibes the local
behavior of the covering spectrum of a manifold
with $sect \ge -1$.  We describe such spaces as having many
``subscaled souls".

In Section~\ref{sect6.3}
 we turn to complete manifolds with nonnegative Ricci curvature.
Theorem~\ref{Riccicutoff} states that the cut-off covering spectrum
of such a space is empty unless its universal cover splits isometrically.
In particular a manifold with positive Ricci curvature has an empty
cut-off covering spectrum [Corollary~\ref{Riccicutcor}].

In Section~\ref{sect6.4} we prove Theorem~\ref{Riccicutoff2}
which concerns limits of spaces with lower bounds on their
Ricci curvature approaching $0$.  We then suggest some
open problems related to Ricci curvature and the cut-off
covering spectrum and possible local almost soul theorems similar to
Corollary~\ref{corosoul2}.  Conjecture~\ref{Conj1} suggests
an extension of a theorem of the first author from \cite{So-loops}
which was used to prove Theorem~\ref{Riccicutoff}.
In Conjecture~\ref{conj2} we suggest
that Theorem~\ref{Riccicutoff2} might then be strengthened to
Theorem~\ref{Riccicutoff}.  In Conjecture~\ref{conj3} we suggest
a possible subscaled soul theorem for manifolds with $Ricci \ge -1$
similar to Corollary~\ref{sectge-epsi}.  We close
by applying Theorem~\ref{covcutconv}
and Theorem~\ref{almostsoul} to prove Conjecture~\ref{conj2}
implies Conjecture~\ref{conj3} [Theorem~\ref{conj2to3}].

Appendix A provides the background on a concept we call {\em semiclosure}
which is needed to describe some properties of the covering spectrum.

In Appendix B we prove Lemmas~\ref{balltoL} and ~\ref{balltoL1} 
which correctly restate and circumvent Lemma 5.8 of \cite{SoWei3}.
While the original Lemma 5.8
was applied to prove Theorem 5.7 of \cite{SoWei3} that
the marked length spectrum of a compact length space with a universal
determines the covering spectrum, these lemmas apply to complete 
length spaces which are not compact and do not have universal covers.
They are used throughout this paper.  The Appendix closes with an
explanation as to why Theorem 5.7 of \cite{SoWei3} is correct.

We would like to thank Carolyn Gordon, David Fisher and 
Ruth Gornet for encouraging us to pursue a further investigation 
of the covering spectrum.  We would like to thank Conrad Plaut
for his incisive questions regarding the existence of rectifiable  
curves that led to Appendix B.  Finally we would like to thank             
Jay Wilkins and the referee for their close reading of the preprint.       
The first author is grateful to the Courant Institute for its             
hospitality in Spring 2007.

\sect{Background} \label{review}\label{sect2}

In Section~\ref{sect2.0} we review the classical notions of
covering spaces and length spaces and how the metric structures
lift from the base space to the covering space.
In Sections~\ref{sect2.1} and~\ref{sect2.2} we review the definitions
of Spanier covers, delta covers and the covering spectrum for
complete length spaces (geodesic spaces).  In 
Section~\ref{sect2.3} we provide a simplified yet
equivalent definition of the 
covering spectrum of a space when the space has a universal cover 
[Defn~\ref{groupcovspec} and Thm~\ref{balltoloop}].  Elements
required for the proof of Theorem~\ref{balltoloop} appear in
Appendix B.

\subsection{Covering Spaces and Length Spaces} \label{sect2.0}

This section reviews classical theory of covering spaces and
length spaces.  First we recall some basic definitions.

\begin{defn}
A metric space is a length space if the distance between
points is the infimum of the lengths of curves running
between those points.  When the infimum is achieved
between a given pair of points, we call the minimizing
curve a minimal geodesic.  When the infimum is achieved
for any pair of points, we say the space is a geodesic space.
\end{defn}

Given a subset $A$ of a length space $X$, one may either use the
restricted metric (which is not a length metric on $A$)
or the {\em induced length metric} which is found by taking the infimum
of lengths curves lying in $A$.  Note that even when $X$ is a geodesic 
space, a subset, $A\subset X$, with the induced length metric might 
not be a geodesic space:

\begin{example}  \label{Aexample}
Let $X$  be the Euclidean plane and 
\be
A=\{(x,y): \, y= |1-x|/j, \, x \in [-1,1], \, j \in \Bbb{N}\}\subset X.
\ee
Then if $d_A$ is the induced length metric,
$d_A((-1,0), (1,0))=2$ is not achieved.
\end{example}

Recall that a metric space is locally compact iff 
every point has a precompact neighborhood.

\begin{theo} [Hopf-Rinow] (c.f. \cite{BBI}\cite{Gr})
\newline
If $X$ is a complete locally compact length space, then
closed balls, $B(x,R)$, in $X$ are compact and 
$X$ is a geodesic space.
\end{theo}

Thus, in particular, complete Riemannian manifolds are geodesic spaces.

Like geodesics in Riemannian manifolds, geodesics in length
spaces are defined to be locally minimizing curves.  Closed geodesics are
geodesics from $S^1$ to the space.  However, geodesics in
length spaces are not
necessarily extendable and they may branch.  This can be seen
for example in the case of a closed disk (where geodesics end at
the boundary) and a tree (where geodesics branch at the vertices).

\begin{example}\label{Hawaii}
A simple example of a complete length space that we will use
repeatedly in this section is a collection of circles of various
radii joined at a point, $p$.  The distance between points on single
circle is just the shorter arclength between them.  Distances
between points, $q_1,q_2$ on distinct circles is the sum of the
shorter arclength from $q_1$ to $p$ and the shorter arclength from
$q_2$ to $p$.  This space is called the Hawaii ring when the
collection of radii is $\{1/j: \, j\in \Bbb{N}\}$.
\end{example}

The following classical definitions can be found for example
in \cite{Sp}.

 \begin{defn} {\em
 We say $\bar{X}$ is a  {\em covering space} of $X$ if there
 is a continuous map $\pi: \bar{X} \to X$ such that $\forall x
 \in X$ there is an open neighborhood $U$ such that
 $\pi^{-1}(U)$ is a disjoint union of open subsets of
 $\bar{X}$ each of which is mapped homeomorphically onto $U$
 by $\pi$ (we say $U$ is evenly covered by $\pi$).   }
 \end{defn}

This is clearly a topological definition which preserves the local
topology.  Thus if $X$ is locally compact, then any cover,
$\bar{X}$, is locally compact as well.  Naturally compactness is
a global condition and does not necessarily lift 
to the cover as can be seen with the
classic example of a line covering a circle.  Nevertheless the covering
space of a compact space must be locally compact.

\begin{defn}\label{univcov} {\em \cite[pp
 62,83]{Sp} We say $\tilde{X}$ is a
 {\em universal cover} of $X$ if $\tilde{X}$ is a cover of $X$
 such that for any other cover $\bar{X}$ of $X$, there is a
 commutative triangle formed by a continuous map $f:\tilde{X}
 \to \bar{X}$ and the two covering projections.}
 \end{defn}

Note that the Hawaii Ring does not have a universal cover.
In \cite{SoWei1}, the authors developed a method of detecting
when a compact space like the Hawaii ring has a universal
cover or not.

If $X$ is a length space, then naturally one can define the lengths
of curves, $C$, on its covering space, $\bar{X}$, by measuring
the lengths of their projections down to $X$: $L(C):=L(\pi\circ C)$.
This defines a metric on $\bar{X}$ which is a length metric, and
it is the unique metric on $\bar{X}$ for which $\pi:\bar{X}\to X$
is a local isometry.  {\em Throughout this paper, we will always use this
metric on the covering spaces of a length space.}

When $X$ is a complete locally compact length space, then $\bar{X}$
is also a complete locally compact length space and, by Hopf-Rinow,
it is a geodesic space.  However, if $X$ is only a geodesic space,
the cover might only be a length space as can be seen in the following
example:  

\begin{example}
Let $X$ be a collection of circles of circumference $2\sqrt{1+(1/j)^2}$
joined at a point.  Let $\bar{X}$ be defined as a collection of
the sets $A$ defined in Example~\ref{Aexample} joined to each other
in a row:
\be
\bar{X}=\{(x+2k,y): \,\, (x,y) \in A, \, k \in \Bbb{Z}\}.
\ee
Then $\pi: \bar{X} \to X=\bar{X}/\Bbb{Z}$ is a regular cover, and
$X$ is a geodesic space, but $\bar{X}$ is only a length space.
\end{example}

Most definitions in this paper are well defined for
arbitrary metric spaces.  However, on occasion we need
a geodesic length structure on the covering spaces or the compactness
of closed balls in the covering spaces, and in those
settings we require that our spaces be complete locally
compact length spaces.  We will also employ local compactness
when applying the Arzela-Ascoli Theorem or the Gromov Compactness
Theorem.

\subsection{Spanier covers and $\delta$-covers} \label{sect2.1}

We now introduce a special collection of covers we will
call Spanier covers as they are described in \cite[Page 81]{Sp}.

\begin{defn}\label{SpanierCover}\label{defspanier}
 Let ${\mathcal U}$ be any collection of open sets covering $X$. For any $p \in
X$, by
 \cite[Page 81]{Sp}, there is a covering space, $\tilde{X}_{\mathcal U}$,
 of $X$ with covering group  $\pi_1(X,{\mathcal U}, p)$, where
 $\pi_1(X,{\mathcal U}, p)$ is a normal subgroup of $\pi_1(X, p)$, generated
 by homotopy classes of closed paths having a representative of the form
 $\alpha^{-1} \circ \beta \circ \alpha$, where $\beta$ is a closed path
 lying in some element of $\mathcal U$ and $\alpha$ is a path from $p$ to
 $\beta(0)$.
\end{defn}

 It is easy to see that a Spanier cover is a regular or Galois cover.
 That is, the lift of any closed loop in Y is either always closed
 or always open in a Spanier cover.  In particular Spanier covers
of a collection of circles of various radii will leave some or none
of the circles as circles and unravel the other circles completely
into a tree.

The following lemma is in Spanier \cite[Ch.2, Sec.5, 8]{Sp}:

\begin{lemma}\label{openrefine}
Let $\mathcal U$ and $\mathcal W$ both be collections of open sets that
cover $X$. Suppose $\mathcal U$ refines $\mathcal W$ in the sense that
for any open set $W$ in $\mathcal W$ there is an open set $U \in
\mathcal U$, such that $U\subset W$. Then the Spanier cover
$\tilde{X}_{\mathcal U}$ covers $\tilde{X}_{\mathcal W}$.
\end{lemma}

Spanier covers will be used to define various covering spaces in
this paper as well as the $\delta$ covers first introduced by the
authors in \cite{SoWei1}.

\begin{defn} \label{defdeltacover}
{\em
 Given $\delta > 0$, the  $\delta$-cover, denoted $\tilde{X}^\delta$,
 of a metric space $X$, is defined to be the Spanier
cover, $\tilde{X}_{{\mathcal U}_{\delta}}$,
 where ${\mathcal U}_\delta$ is the open covering of $X$ consisting of
 all open balls of radius $\delta$.

 The covering group will be denoted $\pi_1(X,\delta, p)\subset \pi_1(X,p)$.
 This is the normal subgroup of $\pi_1(X, p)$, generated
 by homotopy classes of closed paths having a representative of the form
 $\alpha^{-1} \circ \beta \circ \alpha$, where $\beta$ is a closed path
 lying in some ball of radius $\delta$ and $\alpha$ is a path from $p$ to
 $\beta(0)$.

}
\end{defn}

In Example~\ref{Hawaii},
the $\delta$ cover of the space consisting of circles of
various sizes glued at a common point, is a covering space which unravels
all the circles of circumference $2\pi r \ge  2\delta$ and keeps the
smaller circles wrapped as circles.  In particular, when $X$ is
the figure eight created by joining one circle of circumference $2\pi$ and
one circle of circumference $4\pi$ at a common point: then
$\tilde{X}^\delta$ is $X$ itself when $\delta > 2\pi$,
it is a real line with circles of circumference $2\pi$ glued at the points
$\{2j\pi :j\in \Bbb{Z}\}$
when $\delta \in (\pi,2\pi]$ and it is the universal cover $\tilde{X}$
when $\delta \le \pi$.

Note that like all covering spaces,
the $\delta$ covers of complete locally compact length spaces
are also complete locally compact length spaces and are thus geodesic
length spaces.  Example~\ref{nolift} is a complete geodesic length
space with a $\delta$ cover that is not an geodesic length space, so
the local compactness condition is required to lift this property
even to $\delta$ covers.

The $\delta$-covers of compact spaces are surveyed quickly in the
background section of \cite{SoWei3}.  There we proved that
$\delta$ covers of metric spaces are monotone in the sense that if
$\delta_1<\delta_2$ then $\tilde{X}^{\delta_1}$ covers
$\tilde{X}^{\delta_2}$ which just follows from
Lemma~\ref{openrefine}. See \cite[Lemma 2.6]{SoWei3}.

If one has a space where balls
of radius $\delta_1$ and $\delta_2$ have the same topology, the
covering spaces are the same.  In fact, for compact spaces, we
proved the $\delta$ covers  are lower semicontinuous in the sense
that for any $\delta_1>0$ there is a $\delta_2<\delta_1$
sufficiently close to $\delta_1$ such that the two delta covers agree
\cite[Lemma 2.7]{SoWei3}.  This is not true for complete noncompact
spaces. In fact, the space of circles of circumference
$2\pi r_i$ joined at a
point have distinct delta covers for each $\delta_i=\pi
r_i$ so that lower semicontinuity fails when there is a sequence
$r_i$ increasing to $r_0$ \cite[Example 2.8]{SoWei3}.

\subsection{Review of the Covering Spectrum} \label{sect2.2}

In \cite{SoWei3} we introduced the covering spectrum
on compact metric spaces which is well defined
on complete noncompact spaces as well.

\begin{defn}  \label{coveringspectrum}
 Given a metric space $X$, the covering spectrum of
 $X$, denoted  CovSpec$(X)$ is the set of all $\delta > 0$ such
 that
 \be
 \tilde{X}^\delta \neq \tilde{X}^{\delta'}
 \ee
 for all $\delta'> \delta$.
  \end{defn}

The covering spectrum of a finite collection of circles of
circumference $2\pi r_i$ joined (glued) at
a common point is $\{\pi r_i\}$.

For a compact length space the covering spectrum is discrete and the
only accumulation point of the covering spectrum that can occur
outside of the covering spectrum is $0$ \cite[Prop. 3.2]{SoWei3}.
This happens for example with the Hawaii Ring where the circles have
circumference $2\pi r_j=2\pi/j$.

The covering spectra of complete noncompact spaces need not be discrete:

\begin{example} \label{noncompat-cov}
The covering spectrum of a complete noncompact length space can be
$(0,\infty)$ as can be seen by joining the uncountable collection of
circles of circumference $2\pi r$ for every $r\in(0,\infty)$ at a common
point.
This same covering spectrum can be achieved by taking a joined
countable collection of circles of circumference $2\pi r$ for every $r\in
\Bbb{Q}$.
\end{example}

The following lemma is a simple exercise on the definition:

\begin{lemma} \label{dec-closure}
If $\delta_j \in CovSpec(X)$ and $\delta_j$ decrease to a positive
limit $\delta_0>0$, then $\delta_0\in CovSpec(X)$.
\end{lemma}

\begin{example}  \label{closure-ex}
Thus the covering spectrum of the joined collection of circles of
circumference $2\pi r_j=2\pi+2\pi/j$, is
$\{\pi(1+1/j): j\in \Bbb{N}\}\cup \{\pi\}$.
In contrast the covering spectrum of the joined collection of circles
of circumference $2\pi r_j=2\pi-2\pi/j$ is just $\{\pi(1-1/j): j\in
\Bbb{N}\}$.
\end{example}

This is just an indication of the complexity one encounters when
studying the covering spectra of complete noncompact spaces. In the
next section we explore this situation, and in subsequence sections
we introduce alternative spectra which detect properties that the
covering spectrum cannot detect on a complete noncompact space.
Further review of the covering spectra of compact spaces will appear
below.

\subsection{The Covering Spectrum and Deck Transforms}\label{sect2.3}

In our prior papers, we did not like to assume the space had a
universal cover in part because we were applying $\delta$-covers and
the covering spectrum to prove the existence of universal covers.
However, if one does assume the existence of a universal cover, then  there
is a fairly beautiful new perspective on the meaning of the covering
spectrum using its relationship with the group of deck
transforms ${\pi}_1(X)$ on the universal cover, $\tilde{X}$.
When the universal cover is simply connected, this group
of deck transforms is isometric to the fundamental group $\pi_1(X,p)$.

Recall that a $\delta$ cover, $\tilde{X}^\delta$, 
is defined using a covering group. $\pi_1(X,\delta,p)$, 
so with a universal cover we have:
\be \label{groupcover}
\tilde{X}^\delta= \tilde{X}/ \pi_1(X, \delta),
\ee
where $\pi_1(X,\delta)\subset \pi_1(X)$.
This provides us with an equivalent definition
for the covering spectrum:

\begin{defn}\label{groupcovspec}
 Given a metric space $X$, with a universal
cover, $\tilde{X}$, the covering spectrum of $X$
is the set of all $\delta > 0$ such
 that
 \be
 \pi_1(X,\delta) \neq \pi_1(X,\delta') \qquad \forall \delta'>\delta
 \ee
when viewed as subsets of $\pi_1(X)$.
\end{defn}

We can now use the existence of the universal cover to
simplify the definition of $\pi_1(X,\delta)$ so that
Definition~\ref{groupcovspec} can be used to quickly
to recover the covering spectrum of a space whose
deck transforms are well understood.  We begin with
the following standard definition.

\begin{defn} \label{DefL(g)}
Given a complete locally compact length space $X$ with universal cover
$\tilde{X}$, for each element $g \in \pi_1(X)$,
its length $L(g)$ is
\be \label{defLeqn}
L(g) = \inf_{\tilde{x} \in \tilde{M}} d (\tilde{x}, g\tilde{x}).
\ee
\end{defn}

\begin{theo} \label{balltoloop}
Given a complete locally compact length space $X$ with a  universal cover,
the $\delta$ covering group $\pi_1(X,\delta)$
is the subgroup of $\pi_1(X)$ generated by
deck transforms $g$ with $L(g) <2\delta$.  
\end{theo}

To prove this theorem we will apply Lemma~\ref{balltoL1}
which states that {\em when $C$ is a rectifiable
loop lying in a ball of radius $\delta$ then $C$ is freely homotopic
to a product of curves of length $< 2\delta$ }.  This
slight restatement of Lemma 5.8 of \cite{SoWei3}
is located in Appendix B.

We also need the following lemma regarding
the rectifiability of loops representing elements of
$\pi_1(X,\delta)$ for spaces with universal covers.  It
is important to note that in the Hawaii Ring,
which does not have a universal cover, there are elements
of $\pi_1(X,\delta)$ which are not rectifiable.  

\begin{lemma}\label{betarect}
Given a complete locally compact length space $Y$ with a  universal cover,
if $g \in \pi_1(X,\delta)$ has positive length, then it
is generated by elements of $\pi_1(X, \delta)$
whose representatives have the form $\alpha^{-1} \circ \beta\circ\alpha$
where $\beta$ are rectifiable and lie in balls of radius
$\delta$.
\end{lemma}

It is important to note that we do not claim that the fundamental
group always has a rectifiable representative.  Recall that spaces
exist which do not have simply connected universal covers \cite{Sp}.
These spaces can have elements in their fundamental group which
do not have rectifiable curves representing them, but we show
that their images in the deck transform group will have rectifiable
representative curves.

\noindent {\bf{ Proof of Lemma~\ref{betarect}:}}
Suppose that $L(g)=L_0>0$.  So for all $x\in \tilde{X}$,
$d(gx,x) \ge L_0$.  We know $g$ is a product of elements
with representatives of the form 
$\alpha^{-1} \circ\beta'\circ \alpha$ where the $\beta'$
might not be rectifiable but are contained in balls
of radius $< \delta$.  

We need only replace each $\beta'$ with a rectifiable 
curve $\beta$ that has the same end points when lifted to the
universal cover and which also fits in the $\delta$ ball.

Note that since $[0,1]$ is compact, the image of $\beta '$ is 
a compact set.  So there
exists an $\epsilon>0$ sufficiently small that 
for any $t\in S^1$, $B_{\beta '(t)}(5\epsilon)\subset B_q(\delta)$.

Now we lift $\beta'$ to the universal cover $\tilde{X}$ and 
partition the lift $\tilde{\beta}'$ into $0=t_0<t_1<t_2<...<t_N=1$
so that each segment $\tilde{\beta}'([t_i,t_{i+1}])$ lies 
in a ball of radius $\epsilon$.  We create a piecewise geodesic, $\tilde{\beta}$ 
in the universal cover which joins these endpoints with minimal
geodesics of length $<2\epsilon$.

Projecting $\tilde{\beta}$ down to a curve $\beta$, we see
$\beta$ has the same endpoints at $\beta'$, passes through
$\beta'(t_i)$ and must also lie in 
$$
\bigcup B_{\beta'(t_i)}(2\epsilon)\subset B_q(\delta).
$$
\qed

We now easily prove the theorem:

\noindent{ \bf{ Proof of Theorem~\ref{balltoloop}:}}
Suppose $L(g)<2\delta$.  Then $g$ has a representative loop
$\beta$ whose lift runs between some $x$ and $gx$ of
length less than $2\delta$. We do not need to assume compactness of $X$
because of the strict inequality here.  This implies $\beta$
has length $< 2\delta$ and thus fits in some ball of radius $\delta$.
So $g \in \pi_1(X,\delta)$.

Suppose on the other hand $g\in \pi_1(X, \delta)$.  
By Lemma~\ref{betarect},  either
$L(g)=0< 2\delta$ or it 
is generated by elements of $\pi_1(X,\delta)$ with representative
loops of the form $\alpha^{-1} \circ \beta \circ \alpha$
where $\beta$ is rectifiable and lies in a ball of radius $\delta$.  
We just apply Lemma~\ref{balltoL1} to state that 
$\beta$ is homotopic to a product of loops of length $< 2\delta$.
Thus $g$ is generated by  elements with representative
loops of length $< 2\delta$.  Such elements must then have length
$< 2\delta$. 
\qed

On compact length spaces 
Theorem~\ref{balltoloop} can be combined with
Arzela-Ascoli and Lemma~\ref{achieve} to prove
$CovSpec(M)\subset(1/2)Length(M)$ where $Length(M)$ is the collection
of lengths of closed geodesics $\gamma:S^1 \to M$ \cite{SoWei3}.  This is
not true on complete locally compact length spaces as the 
infimum in (\ref{defLeqn}) need
not be achieved:

\begin{example}\label{notslip}
Let
$M^2$ be the warped product manifold $\Bbb{R} \times_{f(r)} S^1$
where
\be
f(r)=2Arctan(-r)+2\pi.
\ee
Here $\pi_1(M)$ is generated by a single element $g$ whose
length
\be
L(g)=\inf_{r\in (-\infty,\infty)} f(r)=\pi,
\ee
but there is no closed curve homotopic to a representative
of $g$ whose length is $\pi$.
\end{example}

On a compact Riemannian manifold,  $CovSpec(M)=\emptyset$ implies
$M$ is simply connected \cite{SoWei3}.  Yet this is not
true for complete manifolds:

\begin{example} \label{slipex}
Let
$M^2$ be the warped product manifold $\Bbb{R} \times_{f(r)} S^1$
where
\be
f(r)=2Arctan(-r)+\pi.
\ee
Given any $\delta>0$, eventually $f(r)<2\delta$, so
any $g\in \pi_1(M)$ is represented by a loop of length $<2\delta$.
Thus by Theorem~\ref{balltoloop}
the covering spectrum is empty.
\end{example}

Further implications of this perspective on the covering
spectrum will be investigated in future joint work.
In that paper we will also
investigate the slipping group:

\begin{defn}
The slipping group of $X$, denoted
$\pi_{slip}(X)$, is generated by the elements $g\in \pi_1(X)$
such that $L(g)=0$.
\end{defn}

\sect{Delta homotopies} \label{delta-homotopy} \label{sect3}

In this section we develop the concept of the delta homotopy which
we first defined in \cite{SoWei1}:

\begin{defn}\label{defdelhom}
Two loops $\gamma_1, \gamma_2$ in a metric space, $X$,
are called $\delta$-homotopic if
$\pi_\delta ([\gamma_1]) = \pi_\delta ([\gamma_2])$,
where $\pi_\delta: \pi(X) \rightarrow \pi(X)/\pi(X, \delta)$.
In particular $\gamma_1$ is $\delta$ homotopic to a point
if
\be
[\gamma_1] \in \pi(X,\delta)
\ee
which means $\gamma_1$ lifts as a closed loop to $\tilde{X}^\delta$.
\end{defn}

This concept can be used to produce closed geodesics in length spaces:

\begin{lemma} \label{achieve}
If $X$ is a length space and $\gamma: S^1 \to X$ 
has $L(\gamma)\le 2\delta$ but is not $\delta$ homotopic to a point 
then $\gamma$ is a closed geodesic which is
minimizing over any interval of half its length and has length
$2\delta$.
\end{lemma}

This lemma will be applied later when we prove our new spectra are
in the length spectrum.

\Pf
Since $\gamma$ lifts as a closed loop to the length space
$\tilde{X}^\delta$ it does not
fit in a ball of radius $\delta$.  In particular, for any
$t\in S^1_{\delta/\pi}$ we have
\be \label{achieve1}
Im(\gamma) \cap (X \setminus B_{\gamma(t)}(\delta)) \neq \emptyset.
\ee
However $L(\gamma)=2\delta$ so the only point in (\ref{achieve1})
must be $\gamma(t+\delta)$ and $d(\gamma(t+\delta), \gamma(t))$ must
be $\delta$.  Thus $\gamma$ is minimizing on any subinterval of
length $\delta$ including an interval centered at $t=0$.
\qed

The remiander of this section will be dedicated to providing a more
geometric understanding of $\delta$ homotopies on metric spaces.  
We will first relate
$\delta$ homotopies to grids [Section~\ref{sect3.1}] and describe
how to localize $\delta$ homotopies [Section~\ref{sect3.2}].
Finally we prove a few properties of $\delta$ homotopies that are
localized in precompact sets in length spaces [Section~\ref{sect3.3}].

While we apply the results in
this section to study the cut-off covering spectrum, we prove them
first because they apply in a much more general setting and should
prove useful for those interested in other concepts.  Those who are
more interested in the cut-off covering spectrum may jump to
Section~\ref{sect4} and only return to this section before continuing to
Section~\ref{sect5} on Gromov-Hausdorff convergence.  Alternatively
one might skim through this section reading only the statements and
viewing the accompanying diagrams.

\subsection{Using grids to understand $\delta$ homotopies} \label{sect3.1}

Before we can transform our original somewhat algebraic definition
of $\delta$ homotopy [Definition~\ref{defdelhom}] into a geometric
statemant about grids on metric spaces
, we need to examine the definition closely.
Clearly it is base point independent.  So if a curve
$C$ is $\delta$ homotopic to a point then $\alpha C \alpha^{-1}$ is
also $\delta$ homotopic to a point.  So it is often easier to think
of $\gamma_1$ as $\delta$ homotopic to $\gamma_2$ if we joint them
to a common point via curves $\alpha_1$ and $\alpha_2$ and then say
$\alpha_1\gamma_1\alpha_1^{-1}$ is $\delta$ homotopic to
$\alpha_2\gamma_2\alpha_2^{-1}$ which is the same as saying \be
\alpha_1\gamma_1\alpha_1^{-1}(\alpha_2\gamma_2\alpha_2^{-1})^{-1}
=\alpha_1\gamma_1\alpha_1^{-1}\alpha_2^{-1}\gamma_2^{-1}\alpha_2 \ee
is $\delta$ homotopic to a point.  In this sense we make the
following definition:

\begin{defn}
A collection of loops $\gamma_1, \gamma_2,...\gamma_k$ is $\delta$
homotopic to a point if there exist curves $\alpha_i$ mapping a base
point $p$ to $\gamma_i(0)$ and such that \be
\alpha_1\gamma_1\alpha_1^{-1}\alpha_2\gamma_2\alpha_2^{-1}...
\alpha_k\gamma_k\alpha_k^{-1} \ee is $\delta$ homotopic to a point.
\end{defn}

The ordering of the loops is important in this definition.
If $\gamma_1, \gamma_2$ is $\delta$ homotopic to a point then
$\gamma_1$ is $\delta$ homotopic to $\gamma_2^{-1}$.

\begin{lemma} \label{deltahomotopyviagrid}
\label{equivdefdelhom}\label{lemma3.3}
A loop $C$ of length $L$ in a metric space
is $\delta$ homotopic to a point iff there
is a $\delta$ homotopy
$H:G \to X$
where $G$ is an $N\times M$ grid of unit squares
such that $H(0,y)=C(yL/M)$, $H(x,0)=H(x,M)=H(N,y)=C(0)$
and such that the image under $H$ of each square in
the grid is contained in a ball of radius $\delta$.
\end{lemma}

\begin{figure}[htbp]
\includegraphics[width=4.5in]{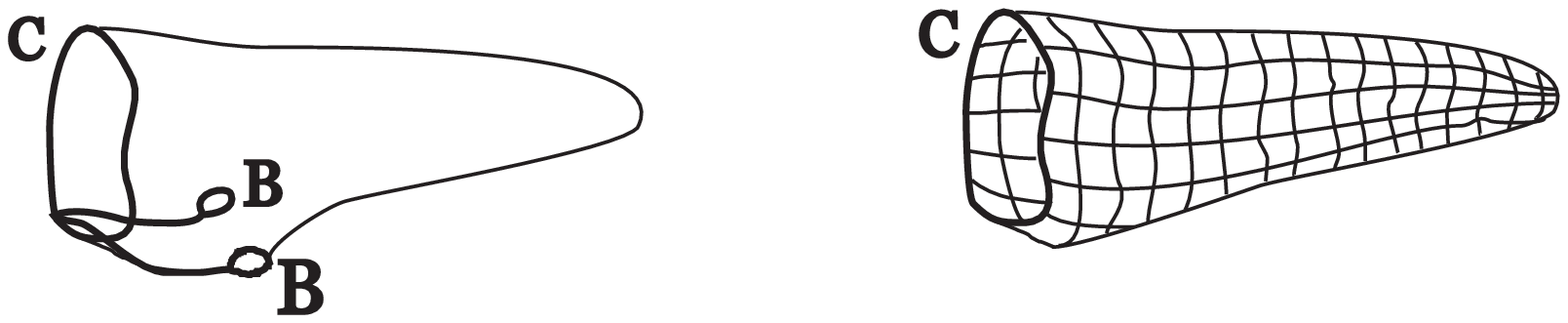}
\caption{}
\label{togrid}
 \end{figure}

In some sense this lemma is intuitively obvious.  See
Figure~\ref{togrid}.  Special cases of this lemma were
used within some of the proofs in \cite{SoWei1}.
Writing out the proof is a bit technical and so
first we set some notation.  Let
$\beta_{j,k}$ be image of the clockwise loop around the square
$(j,k),(j,k+1),(j+1,k+1), (j+1,k)$.  Let $\alpha_{j,k}$
be the image of the line segment from $(j,0)$ to $(j,k)$.
Let $\bar{\alpha}_j$ is the image of the line segment
from $(j,0)$ to $(j-1,0)$.

\Pf If such a homotopy exists, then define $C_j(t)$ to be the loop
$H(j,t)$ from $t$ to $M$ so $C_0(t)=C(tL/M)$ and $C_N(t)$ is a
point. Note that $C_0$ is just $C$. Furthermore each \be
\alpha_{j,0}\beta_{j,0}\alpha_{j,0}^{-1}
\alpha_{j,1}\beta_{j,1}\alpha_{j,1}^{-1}...
\alpha_{j,M}\beta_{j,M}\alpha_{j,M}^{-1} \textrm{ is homotopic to }
C_j(t)(\bar{\alpha}_j C_{j-1}(t) \bar{\alpha}_j^{-1})^{-1} \ee
within the image of the grid.  Thus by the definition of the
$\delta$ cover, \be C_j(t)(\bar{\alpha}_j C_{j-1}(t)
\bar{\alpha}_j^{-1})^{-1} \ee lifts as a closed loop to the
$\tilde{X}^\delta$ and so $C_j$ and $C_{j-1}$ are $\delta$ homotopic
to each other. Thus $C$ is $\delta$ homotopic to $C_N$ which is a
point.

Conversely, if $C$ is $\delta$ homotopic to a point, then
by the definition of the $\delta$ cover, $C$ is homotopic
to a collection of curves $\alpha_i\beta_i\alpha_i^{-1}$
where $\beta_i$ are in balls of radius $\delta$.  So we
take the homotopy $\bar{H}:[0,N]\times[0,M]\to X$
so that $\bar{H}(0,t)=C(tM/L)$,
$\bar{H}(s,0)=\bar{H}(s,1)=C(0)$ and
\be
\bar{H}(N,t) = \alpha_1\beta_1\alpha_1^{-1}\alpha_2\beta_2\alpha_2^{-1}...
               \alpha_k\beta_k\alpha_k^{-1}(t)
\ee
Using the uniform continuity of the homotopy $\bar{H}$
we can choose $N$ and $M$ large enough that each square
in the grid is within a ball of radius $\delta$.  We can
also insure, possibly by adding a few more columns to
allow for a slow homotopy between reparametrizations,
that each $\beta_j$ starts at a $t_j$ and ends
at a $t_j+1$ where $t_j$ are integers.

We now add a gridded column of unit squares on the
right side of the homotopy.  The horizontal bars
will have constant images.  The verticals will agree
with $\bar{H}(N,t)$ whenever this is part of an $\alpha$
curve but will take the value $\bar{H}(N,t_j)$ for
$t\in [t_j, t_j+1]$.  In this way most of the new squares
will be in subsegments of the $\alpha$ curves, and the
selected new squares at the $t_j$ points will have images
equal to $\beta_j$ and thus lie in balls of radius
$\delta$.

Finally we add a number more columns to allow for a homotopy
from the curve
\be
\alpha_1\alpha_1^{-1}\alpha_2\alpha_2^{-1} ...\alpha_k\alpha_k^{-1}
\ee
to a point.  This can be done just by contracting along each
$\alpha_j$.  In this way we complete the homotopy.
Then we restrict the homotopy to the grid points and we
are finished.
\qed

\begin{lemma} \label{squareinball}
If $H$ is a $\delta$ homotopy to a metric space
, then there exists $\epsilon \in (0,\delta)$
sufficiently close to $\delta$ that $H$ is an $\epsilon$ homotopy.
\end{lemma}

In fact on compact spaces, one then has
$\tilde{X}^\epsilon=\tilde{X}^\delta$ as proven in Lemma 2.7 of
\cite{SoWei3}.

\Pf
By Definition~\ref{extdefdelhom}, every square $S_{i,j}$ in the domain of
$H$
is mapped into a ball $B_{q_{i,j}}(\delta)$.  Since $H(S_{i,j})$ is
a closed set, lying in an open ball, it fits in a smaller open
ball $B_{q_{i,j}}(\delta_{i,j})$ with $\delta_{i,j}<\delta$.
Let
\be
\epsilon=\max \{\delta_{i,j}:\, i=1..N, \, j=1..M\}.
\ee
\qed

\subsection{$\delta$ homotopies in subsets} \label{sect3.2}

The following extension of the definition of $\delta$ homotopy
takes full advantage of Lemma~\ref{deltahomotopyviagrid}.  Note
that this extension only requires $X$ to be a metric space.

\begin{defn} \label{extdefdelhom}
A loop $C$ of length $L$ in a metric space, $X$,
is $\delta$ homotopic in $A \subset X$  to a point
if there
is a $\delta$ homotopy
$H:G \to A$
where $G$ is an $N\times M$ grid of unit squares
such that $H(0,y)=C(yL/M)$, $H(x,0)=H(x,M)=H(N,y)=C(0)$
and such that the image under $H$ of each square in
the grid is contained in a ball of radius $\delta$.

We say a curve $C_0$ is $\delta$ homotopic in $A$
to a collection of loops $C_1,C_2,...,C_k$, if
there exists paths $\alpha_j$ from
$C_0(0)$ to $C_j(0)$ lying in $A$ such that
\be
C_0^{-1}\alpha_1C_1\alpha_1^{-1} \alpha_2
C_2\alpha_2^{-1} ...\alpha_kC_k\alpha_k^{-1}
\ee
is $\delta$ homotopic in $A$ to a point. We will say $C_0$ is $\delta$
homotopic in $A$ to a collection of loops in $B \subset A$ if the loops
$C_j$ lie in $B$ but we do not require the
paths $\alpha_j$ lie in $B$.
Similarly one can define $\delta$ homotopies in $A$
between two collections of curves.
\end{defn}

Suppose we have a curve which is $\delta$ homotopic
in a set $A$ to a point and we would like to restrict the
$\delta$ homotopy to a set $B \subset A$.  Parts of
the $\delta$ homotopy may well leave $B$ and so they
need to be chopped off.  This provides new curves where
the homotopy is chopped.  See Figure~\ref{cutnet}
for a glimpse of an application.

\begin{lemma} \label{squaregridlem} \label{lemma3.5}
Given a $\delta$ homotopy in $A$, $H:G \to A$ from a
curve $\gamma$ to a point, and given a set $B$ contained
in $A$ such that $\gamma \subset B$, then
$\gamma$ is $\delta$ homotopic in $B$ to a collection
of curves $\gamma_1, \gamma_2, \cdots ,\gamma_k$ such that
each $\gamma_j$ lies in
$ B$ and  the tubular neighborhood $T_{2\delta}(A\setminus B)$.
\end{lemma}


\begin{figure}[htbp]
\includegraphics[width=5in ]{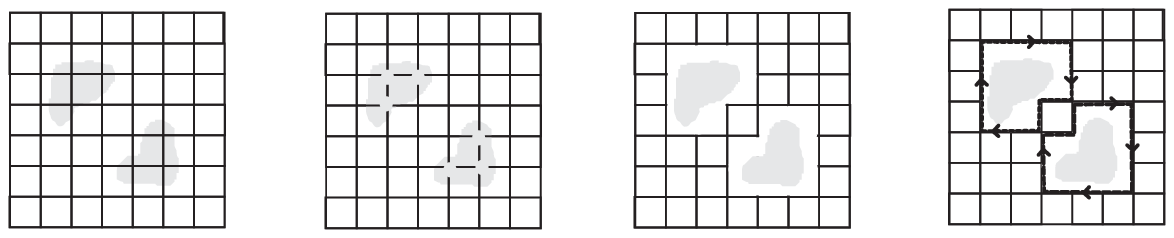}
\caption{}
\label{figsquaregrid}
 \end{figure}

Figure~\ref{figsquaregrid} depicts this lemma and the idea of the
proof. The grey regions are the pullback of the $A\setminus B$ to
the grid of the initial homotopy.  In the figure the collection is
just a pair of curves.  Once we have the final picture in the
figure, we can apply Lemma~\ref{deltahomlem} stated and proven below
to justify that images of the pair of curves
produced in the last step of the picture are indeed $\delta$
homotopic to the initial curve.

\begin{lemma} \label{deltahomlem}
Given a $\delta$ homotopy in $A$, $H:G \to A$, from a curve $\gamma$ to a
point, and a subset of squares $G'\subset G$ such that
the image of $Cl(G \setminus G')$ is contained in a set $B\subset A$.
Here by closure, we are including the boundary of $G'$.
Suppose $G'$ has connected components $G_1,...G_k$.  Let
$\gamma_j$ be the boundary of $G_j$ running around clockwise
so that the image of $\gamma_j$ lies in $B$.

Then $\gamma$ is $\delta$ homotopic in $B$ to the collection of
curves $\alpha_j\gamma_j \alpha_j^{-1}$ where $\alpha_j$ are
paths lying in $B$ or, equivalently, is
freely $\delta$ homotopic in $B$ to the collection of
curves $\gamma_j$.
\end{lemma}

\begin{figure}[htbp]
\includegraphics[width=4.5in ]{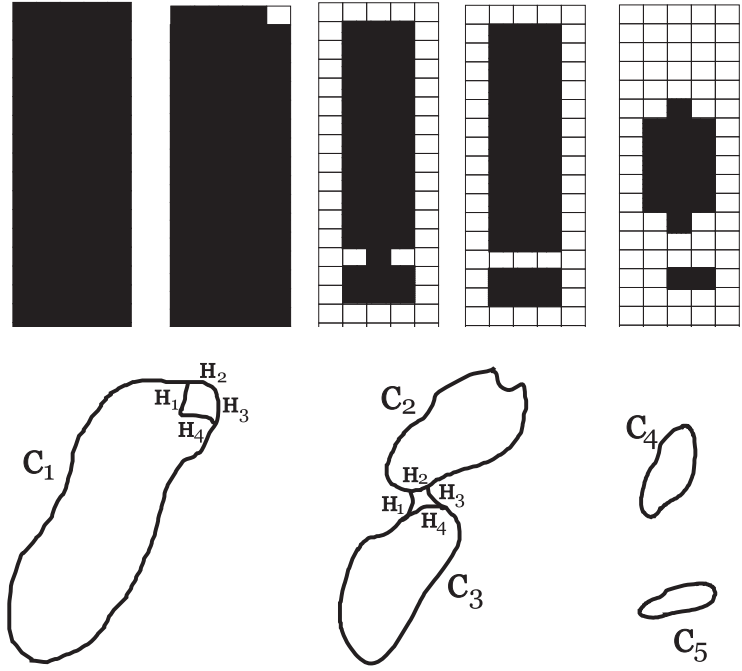}
\caption{In this figure $\gamma_0=C_1H_2H_3$ and we are proving it
is $\delta$ homotopic in $B$ to the pair of curves $\gamma_1=C_4$
and $\gamma_2=C_5$.  The grids are drawn above.  The first rectangle
is filled in completely so we can view our $\gamma_0$ as the
boundary of the full dark region.  The last rectangle has $G'$
darkened and it's two connected components $G_1$ above $G_2$.
Intuitively we are saying that the two inner curves should be
$\delta$ homotopic to the outer curve because of all the squares
between them.  The rectangles in between show how we can run through
a sequence of subsets of the grid creating a $\delta$ homotopy from
$\gamma_0$ to the pair $\gamma_1$ and $\gamma_2$.} \label{deltahom}
 \end{figure}

Intuitively this can be seen because there are only squares that fit in
balls of radius $\delta$ running between them.  You might
wish to skip the proof if you intuitively believe the
process.  For the intuitive idea see Figure~\ref{deltahom}.

\Pf We now rigorously construct a sequence of collections of curves
so that each collection is $\delta$ homotopic to the next.  We begin
with $\gamma_0$ which is the image of the boundary of the entire
grid $G_0=G$.  Each $G_i$ will be a subset of $G_{i-1}$ created by
removing one square, and at each step our collection of curves will
be the boundary of $G_i$.  We know that we can create a sequence of
$G_i$ so that eventually we arrive at $G_I=G'$.  We just need to
verify that we have a $\delta$ homotopy running from each boundary
to the next. There are three cases.

The first case we encounter occurs when
removing a square does not change the number of connected
components of the subgrid.  This is seen in the first part
of Figure~\ref{deltahom}.  A square is removed from the
side on one region.  We need to show that a curve of the
form $C_1H_2H_3$ is $\delta$ homotopic to $C_1H_1^{-1}H_4^{-1}$
when $H_1,H_2,H_3H_4$ is a loop in a ball of radius $\delta$
because it is the image of a single square.  To construct
the $\delta$ homotopy, we set $H(0,t)$ to be the required
\be
C_1H_2H_3(C_1H_1^{-1}H_4^{-1})^{-1}=
C_1H_2H_3H_4H_1C_1^{-1}.
\ee
This time we put all of $H_2H_3H_4H_1$ into one integer
segment and stretch the $C_1$ enough that each segment
lies in a $\delta$ ball.  We add the second column to the
grid keeping everything as in the first column except for
the $H_2H_3H_4H_1$ segment which is now just set to $H_2(0)=C_1(L)$.
Thus the image of the grid thus far is contained in the
images of the old curves which is in $B$ and all the
squares are in $\delta$ balls trivially.  The rest of
the homototy is a classical homotopy contracting
$C_1C_1^{-1}$ to a point and we take as many columns as necessary
so that everything is done slowly enough to fit in
balls of radius $\delta$.  This portion is contained
in $Im(C_1)\subset B$ so we are done.
It is also possible that the square would
be attached on only one side, but this is equally easy.

The second possible case, depicted in the center of
Figure~\ref{deltahom} is when the square which is removed creates
divides a region into two connected components.  So we must show
that $C_3H_1C_2H_3$ is $\delta$ homotopic in $B$ to the pair of
curves $C_3H_4^{-1}$ and $C_2H_2^{-1}$ given that $H_1H_2H_3H_4$ is
the image of a square and so lies in a ball of radius $\delta$.
Using $H_4H_1$ to run $C_2H_2$ to a common base point, we will
construct a $\delta$ homotopy in $B$ from \be
(C_3H_1C_2H_3)((H_4H_1)C_2H_2^{-1}(H_4H_1)^{-1})^{-1}(C_3H_4^{-1})^{-1}
\ee to a point.This is already homotopic within its image to \be
C_3H_1C_2H_3(H_4H_1)H_2C_2^{-1} (H_4H_1)^{-1}H_4C_3^{-1} \ee which
is homotopic within its range to \be C_3H_1C_2H_3(H_4H_1)H_2C_2^{-1}
H_1^{-1}C_3^{-1} \ee Once again we set this up as the first column
so that each collection of curves $H_j$ fit in a single unit segment
and the $C_j$ are spread out so that they divided into pieces of
length less than $\delta$. Our second column will be set up so that
all the horizontal bars are constant and the new vertical line is
the same as before except that the segment with $H_3(H_4H_1)H_2$ is
not just the fixed point $H_3(0)=C_2(L_2)$.  So our new column is
\be C_3H_1C_2C_2^{-1} H_1^{-1}C_3^{-1} \ee but this can be
contracted via a homotopy lying on its image to a point, so we just
provide that homotopy enough columns so that the images of all the
squares lie in $\delta$ balls.  So we are done. Note that the order
of the new collection was important so that this last step would
untangle.

In fact there are cases where the square that is removed might separate into
three or even four regions. This follows exactly as above the the
regions need to be selected in clockwise order around the square to
get the last step to untangle.

The third case is the situation where removing a square removes a
segment from the collection. That situation is trivial.  Anytime a
collection of curves includes a loop within a $\delta$ ball it is
$\delta$ homotopic to the collection with the ball removed.

Thus we have shown that no matter how we remove the
square, we can show that each collection of curves
is $\delta$ homotopic to the next collection carefully
replacing one curve by a new curve or a new curve by
a collection of new curves in the right order until
finally one has the boundary of the given region $G'$.
\qed

We can now return to the proof of Lemma~\ref{squaregridlem}.
See Figure~\ref{figsquaregrid}.

\noindent{\bf Proof of Lemma~\ref{squaregridlem}:}
Let $H$ be the given $\delta$ homotopy.
Remove all vertices in $G$ which are mapped by $H$ into $A \setminus
B$. Remove all the squares touching these vertices.  This gives our
collection of squares $G'$ which satisfies the condition of
Lemma~\ref{deltahomlem}.  So we obtain a collection
$\gamma_1,...\gamma_k$ which are $\delta$ homotopic in $B$ to
$\gamma$ where each $\gamma_j$ lies in the boundary of $G'$.  Thus
every point $q$ which lies on a $\gamma_j$, is on the image of a
square which includes one of the original removed points $z$.  So
$q$ and $z$ lie in a common ball of radius $\delta$ and $z\in A
\setminus B$.  Thus $q\in T_{2\delta}(A\setminus B)$. \qed

\subsection{Compactness and $\delta$ homotopies} \label{sect3.3}

One very nice attribute of $\delta$ homotopy classes of curves is
that they interact well with compactness on length spaces
so that one can control
the lengths of curves in a given class. See Section~\ref{sect2.0}
for a review of these concepts.

\begin{lemma} \label{lengthbound}\label{lem3.7}
Let $\mathcal{C}$ be a set of loops in a precompact length space,
$Z$, which includes a trivial loop.
Suppose there is a curve $C$ in $Z$
which is not $\delta=5\rho$ homotopic to
any collection of curves in $\mathcal{C}$.
Suppose the
number of disjoint balls of radius $\rho$ lying in
$Z$ is bounded above by a finite
number $N$.  Then there exists a curve $\gamma$
in $Z$ which is not
$\delta=5\rho$ homotopic to
any collection of curves in $\mathcal{C}$ and has length
$\le 5N\rho$.
\end{lemma}

Note here that we cannot just take $C$ to be a trivial loop (a
point), because then it would be homotopic to the trivial loop in
$\mathcal{C}$.  In our application $\mathcal{C}$ will be all loops located
outside a given set but this more general statement is equally valid
and possibly useful to others.  Our space $Z$ will be a subset of a
larger space using the induced length metric and thus might not be
complete.

\Pf
Take a maximal disjoint collection of balls of
radius $\rho$
centered at points
\be
Y=\{y_j: j=1..N\}\in Z.
\ee
So the tubular neighborhood of radius $2\rho$
of this finite collection of points contains all of $Z$.

Take $C:[0,L]\to Z\subset T_{2\rho}(Y)$
parametrized by arclength which is not $5\rho$
homotopic to any collection of curves in $\mathcal{C}$.
We will use $C$ to construct a shorter such curve.
Define $0=t_0<t_1< \cdots < t_k=1$ such that
$t_j-t_{j-1}=\rho$ for $j<k$ and
$t_k-t_{k-1}<\rho$.  Define
$\sigma:[0,L]\to Z$ so that $\sigma(t_j)$ is a point
in $Y$ closest to $C(t_j)$.  Then
\be
d_A(\sigma(t_j),\sigma(t_{j+1})\le
d_A(\sigma(t_j),C(t_{j})+\rho
d_A(C(t_{j+1}),\sigma(t_{j+1}) < 5\rho
\ee
and we can join the points in $\sigma$ by curves in $A$
of length $< 5\rho$.  We can also join
$C(t_j)$ to $\sigma(t_j)$ by a curve $h_j$ in $A$ of
length $<2\rho$.  Thus we have a collection of
squares
\be
h_j \sigma([t_j,t_j+1]) h_{j+1}^{-1} C([t_j,t_{j+1}])^{-1}
\subset B_{\sigma(t_j)}(5\rho).
\ee

So $C$ is $5\rho$ homotopic to $\sigma$.
Thus $\sigma$ is not $5\rho$
homotopic to
any collection of curves in $\mathcal{C}$.
If $k\le N$
then $L(\sigma) \le 5k\rho \le 5N\rho$ and
we are done.

If $k>N$ then by the pigeon hole principle
and the fact that
$\sigma(t_j) \in Y$ for $j=0..k$.
with $\sigma(t_0)=\sigma(t_k)$. We see that
there must exist a pair $m,n \in \{0,...k-1\}$
with $|m-n|\le N$ such that
$\sigma(t_n)=\sigma(t_m)$.  This allows us
to break our loops $\sigma$ into two loops
one of which is of length $\le 5\rho N$.
If the other loop is longer, apply the
pigeon hole principle to that loop, and
break off another loop of length $\le 5\rho N$.
Repeating this at most finitely many times,
we see that our original curve $\sigma$ is really
a concatenation of loops all of which have length
$\le 5 \rho N$.

We claim one of these short loops must not be
$5\rho$
homotopic to any collection of curves in $\mathcal{C}$.
Otherwise, all the of them are
$5\rho$
homotopic to some collection of curves in $\mathcal{C}$
and so their concatenation must be $5\rho$
homotopic to a concatenation of that collection.
\qed

\sect{The Cut-off Covering Spectrum} \label{cut-off-spectrum}\label{sect4}

It is natural when studying complete noncompact manifolds
to remove the ends of the manifolds before beginning the
analysis.  This technique works well for any complete locally compact 
metric spaces. In fact, it is standard to refer to pointed
spaces, $(X,x)$, with a special base point $x\in X$.
In this vein of thought, we define the
cut-off covering spectra.  We begin by defining the
$R$ cut-off $\delta$ covers and $R$ cut-off covering spectra,
$CovSpec^R_{cut}(X)$, which
are blind to everything outside a fixed ball of
radius $R$ as trivial.  Next we define the cut-off $\delta$ covers
by taking $R \to \infty$ and define the cut-off covering spectrum,
$CovSpec_{cut}(X)$, based on them.

While the covering spectrum is not well related to the length
spectrum on length spaces which are only complete and locally
compact, as was seen in
Example~\ref{notslip}, we do prove the $CovSpec_{cut}^R(X)\subset
(1/2)L(X)$ and $CovSpec_{cut}(X)\subset Cl_{lower}((1/2)L(X))$
on such spaces.  We
then review the loops to infinity property, and prove such loops are
not detected by the cut-off covering spectra. We close the section
with two technical subsections: one establishing that the $R$
cut-off covering spectrum is truly localized and the other
describing how $CovSpec^R_{cut}(X)$ changes as one varies $R$.
These results will be applied to establish the continuity properties
of these cut-off spectra in Section~\ref{sect5}.

\subsection{The $R$ cut-off $\delta$ covers and
$CovSpec^R_{cut}(X)$}\label{sect4.1}

The $R$ cut-off covering spectrum  is a basepoint dependant concept.
It is defined on pointed metric spaces $(X,x)$ which are metric spaces
with given basepoints.  We begin with the corresponding covering spaces.
Recall Defn~\ref{defspanier}
of a Spanier Cover.

\begin{defn}\label{cutoffcover}
Given a pointed metric space $(X,x)$,
the $R$ cut-off $\delta$ cover based at $x$, denoted
$\tilde{X}^{\delta, R}_{cut}$ or $\tilde{X}^{\delta, R}_{cut \, x}$,
is the Spanier cover
corresponding to the open sets
\be
\{B_p(\delta):p\in X\}\cup \{ X \setminus \bar{B}_{x}(R) \}.
\ee
\end{defn}

When the basepoint is obvious we will omit it.

\begin{lemma} \label{deltatop}
The $R$ cut-off $\delta$ cover based at $x$
is covered by the $\delta$ cover.
In fact
\be
\tilde{X}^{\delta,R}_{cut}=\tilde{X}^\delta/G(R)
\ee
where $G(R)$ is the subgroup of $\pi_1$ generated
by elements with representative loops
of the form $\alpha\circ\beta\circ\alpha^{-1}$
where $\beta\in M\setminus \bar{B}_x(R)$.
\end{lemma}

\Pf
By definition $\pi_1(\tilde{X}^{\delta,R}_{cut})$
is generated by loops of the form
$\alpha\circ\beta\circ\alpha^{-1}$ where
$\beta$ is either in a ball of radius $\delta$
or in $M\setminus \bar{B}_x(R)$.  So
it is generated by elements in $\pi_1(\tilde{X}^\delta)$
and elements in $G(R)$.  Thus
\be
\tilde{X}^{\delta,R}_{cut}=\tilde{X}/\pi_1(\tilde{X}^{\delta,R}_{cut})
=(\tilde{X}/\pi_1(\tilde{X}^{\delta}))/G(R)
=\tilde{X}^\delta/G(R).
\ee

\begin{lemma}\label{ontop}
If $B_{x_1}(R_1) \subset B_{x_2}(R_2)$  in a metric space $X$
and $\delta_1\le \delta_2$,
then
$\tilde{X}^{\delta_1,R_1}_{cut}$ based at $x_1$ covers
$\tilde{X}^{\delta_2,R_2}_{cut}$ based at $x_2$.
\end{lemma}

\Pf Just apply Lemma~\ref{openrefine} which is proven in \cite{Sp}.
\qed

\begin{example}
A cylinder is its own $R$ cut-off $\delta$ cover for
all $R>0$ and all $\delta>0$.
\end{example}

\begin{defn}\label{defRcutoffspec} \label{defn4.4}
Given a pointed metric space $(X,x)$,
the $R$ cut-off $\delta$ spectrum, denoted $CovSpec^R_{cut}(X)$
or $CovSpec^R_{cut}(X,x)$,
is the collection
of $\delta>0$ such that
\be
\tilde{X}^{\delta_1, R}_{cut} \neq \tilde{X}^{\delta, R}_{cut}
\ee
for all $\delta_1>\delta$.
\end{defn}

Note that by Lemma~\ref{ontop} and Theorem~\ref{appendixthm},
$CovSpec^R_{cut}(X)$ is a lower semiclosed 
set for any metric space $X$.

The following lemma was known for compact length spaces in \cite{SoWei3}:

\begin{lemma} \label{deltaRgap}
Given a pointed metric space $X$,
if $[\delta_1,\delta_0) \cap CovSpec^R_{cut}(X)=\emptyset$,
then $\tilde{X}^{\delta_1,R}_{cut}=\tilde{X}^{\delta_0,R}_{cut}$.
\end{lemma}


\Pf
Let
\be
A=\{\delta \in [\delta_1,\delta_0):
\tilde{X}_{cut}^{\delta,R}=\tilde{X}_{cut}^{\delta_1,R}\}\subset
[\delta_1,\delta_0).
\ee
Claim: $\sup \{A\} = \delta_0$. Otherwise $\sup \{A\} = \delta'< \delta_0$.
By
assumption, $\delta'  \not\in CovSpec_{cut}^R(X,x)$.  Therefore there is
$\delta'' >\delta'$ such that
$\tilde{X}_{cut}^{\delta',R}=\tilde{X}_{cut}^{\delta'',R}$, contradicting
that $\delta'$ is the supremum.

So there exist $\delta_i$ increasing to $\delta_0$ such that
\be \label{oneCenough}
\tilde{X}^{\delta_1,R}_{cut}=\tilde{X}^{\delta_i,R}_{cut}.
\ee

To prove the lemma, we proceed by contradition, assuming
\be
\tilde{X}^{\delta_0,R}_{cut} \neq
\tilde{X}^{\delta_1,R}_{cut}.
\ee
Then there is a curve $C$ which lifts closed to
$\tilde{X}^{\delta_0,R}_{cut}$ but
open to $\tilde{X}^{\delta_1,R}_{cut}$.
Then $C$ is $\delta_0$ homotopic to a collection
of curves outside $\bar{B}(x,R)$.
Applying Lemma~\ref{squareinball} we know that
for $\delta_i$ sufficiently close to $\delta_0$,
$H$ is $\delta_i$ homotopy.  So
$C$ lifts closed to
$\tilde{X}^{\delta_i,R}_{cut}$.
By (\ref{oneCenough}), $C$ lifts closed to
$\tilde{X}^{\delta_1,R}_{cut}$
which is a contradiction.
\qed

\subsection{The cut-off $\delta$ covers and
$CovSpec_{cut}(X)$}\label{sect4.2}

We now introduce a cover which will will later prove is
basepoint independant whenever it is well defined:

\begin{defn}
The cut-off $\delta$ cover of $X$,
denoted $\tilde{X}^\delta_{cut}$ is the Gromov-Hausdorff
limit of the $R$ cut-off $\delta$ covers as $R\to \infty$.
\end{defn}

We do not claim that the cut-off $\delta$ cover is defined 
for an arbitrary metric space, but in the next proposition
[Prop~\ref{cutcoverexists}]
we prove they do exist for complete locally compact length
spaces.  We believe they exist for a much larger class
of spaces but will not be investigating this question ourselves.
For a review of Gromov-Hausdorff convergence see Section~\ref{Sect5}.

Note that
as in the case with the cylinder, whose $R$ cut-off $\delta$
covers are all just they cylinder itself, the cutoff
$\delta$ cover is also just the cylinder.  This is in
contrast with the $\delta$ cover which is Euclidean
space for small enough values of $\delta$.

\begin{prop} \label{cutcoverexists} \label{prop4.6}\label{prop4.7}
For any complete locally compact
length space, the Gromov-Hausdorff limit of the $R$
cut-off $\delta$ covers as $R\to \infty$ exists and does not depend
on the base point $x$. Furthermore we have the following covering
maps: \be \label{cce1} \tilde{X}^\delta \mapsto
\tilde{X}^\delta_{cut} \mapsto \tilde{X}^{\delta,R}_{cut}. \ee
\end{prop}

\Pf First we fix a base point $x\in X$. By Lemma~\ref{deltatop} we
have a sequence of covering maps \be f_R: \tilde{X}^\delta \to
\tilde{X}^{\delta,R}_{cut} \ee and a sequence of covering maps \be
h_R: \tilde{X}^{\delta,R}_{cut} \to X \ee both of which are
isometries on balls of radius $\delta$. 

Let the maximal number of disjoint
balls
of radius $\epsilon$ in a ball of radius $r$ in a space $Y$ be
denoted $N(\epsilon,r,Y)$.   As discussed in Section~\ref{sect2.0},
when $X$ is a complete locally compact
length space, so is $\tilde{X}^\delta$ and so closed
balls in $\tilde{X}^\delta$ are compact.  Thus
$N(\epsilon,r,\tilde{X}^{\delta})<\infty$.
By the covering maps, we then have 
\be
N(\epsilon, r,\tilde{X}^{\delta,R}_{cut})\le N(\epsilon,
r,\tilde{X}^{\delta}) \ee 
so by Gromov's Compactness Theorem, a
subsequence $\tilde{X}^{\delta,R_j}_{cut}$ converges.  We call the
limit $\tilde{X}^{\delta}_{cut}$.  

Furthermore, by the Grove-Petersen
Arzela-Ascoli Theorem, subsequences of $f_{R_j}$ and $h_{R_j}$
converge to functions $f$ and $h$ such that
\be
f: \tilde{X}^\delta
\to \tilde{X}^{\delta}_{cut}
\ee
\be h: \tilde{X}^{\delta}_{cut} \to X
\ee
which are still isometries on balls of radius $\delta/2>0$ and
are thus covering maps.   In fact they are regular covers.
 This implies that any limit space
satisfies (\ref{cce1}).

To show we have a unique limit that doesn't depend on the base
point or the sequence $R_j \to \infty$
, take an alternate base point $x'$ and an alternate sequence
$R'_j \to \infty$.  By the above, a subsequence
converges to some other pointed limit space
$(Z',z')$.  If we call our original limit $(Z,z)$, we now
prove $Z$ is isometric to $Z'$.  

Taking a subsequence so that \be B_x(R_j)\subset
B_{x'}(R'_j)\subset B_x(R_{j+1}) \ee and applying Lemma~\ref{ontop}
we have covering maps \be f_j: \tilde{X}^{\delta,R_{j+1}}_{cut} \to
\tilde{X}^{\delta,R'_j}_{cut} \ee \be h_j:
\tilde{X}^{\delta,R'_j}_{cut} \to \tilde{X}^{\delta,R_j}_{cut} \ee
which are isometries on $\delta$ balls.  Subsequences converge by
Grove-Petersen Arzela-Ascoli to covering maps: 
\be f_\infty:
Z \to Z' \textrm{ and } h_\infty: Z' \to Z. 
\ee  
Note further that $h_\infty(f_\infty(z))=z$ and that
$h_\infty \circ f_\infty$ is an isometry on balls of radius
$\delta/3$ so it is a covering map..

If $f_\infty$ and $h_\infty$ do not form an isometry
and its inverse, then we may assume without loss of
generality that $f_\infty$ is not a one-to-one cover and
there exists $w$ such that $f_\infty(z)=f_\infty(w)$.
Let $\gamma_1$ be a minimizing geodesic between $z$
and $w$ of length $L=d_Z(z, w)$.  

Note that by Hopf-Rinow, the closed ball $B(z,L)$
is compact.  Thus we can lift $\gamma_1:[0,L] \to B(z,L)$ 
to a curve $\gamma_2: [0,L] \to B(z,L)$ such that $\gamma_2(0)=z$,
$L(\gamma_2)=L$ and $h_\infty(f_\infty(\gamma_2(t)))=\gamma_1(t)$.
We may repeatedly lift the curves to obtain $\gamma_j: [0,L] \to B(z,L)$
such that $\gamma_j(0)=z$,
$L(\gamma_j)=L$ and $h_\infty(f_\infty(\gamma_j(t)))=\gamma_{j-1}(t)$.
By the compactness of $B(z,L)$ we know that a subsequence of the
points $\gamma_j(L)$ must converge.  In particular, there exists
$k>l \in \Bbb{N}$ such that 
\be \label{jlthis}
d=d_Z(\gamma_k(L), \gamma_l(L))<\delta/3.
\ee

Since $h_\infty\circ f_\infty$ is an isometry on balls of radius
$\delta/3$, we may then apply this map arbitrarily many times and the
images of these two points will still be a distance $d$ apart.
Thus
\be 
d=d_Z((h_\infty\circ f_\infty)^{k-1}\circ\gamma_k(L), 
(h_\infty\circ f_\infty)^{k-1}\circ\gamma_l(L))=d_Z(z,w)
\ee
and
\be 
d=d_Z((h_\infty\circ f_\infty)^{k}\circ\gamma_k(L), 
(h_\infty\circ f_\infty)^{k}\circ\gamma_l(L))=d_Z(z,z)=0.
\ee
Thus $z=w$ and  we have an isometry between $Z$ and $Z'$
and $\tilde{X}^{\delta}_{cut}$ is uniquely defined.
\qed

We leave the following proposition as an exercise as it can be proven
using similar limits of covering maps:

\begin{prop} \label{cutcoveringorder}
For all $\delta_1 <\delta_2$ we have
\be \label{cce2}
\tilde{X}^{\delta_1}_{cut} \mapsto \tilde{X}^{\delta_2}_{cut}.
\ee
\end{prop}

\begin{defn}\label{defcutoffspec} \label{defn4.8}
The cut-off covering spectrum, denoted $CovSpec_{cut}(X)$,
is the collection of $\delta>0$ such that
\be
\tilde{X}^{\delta_1}_{cut} \neq \tilde{X}^{\delta}_{cut}
\ee
for all $\delta_1>\delta$.
\end{defn}

This spectrum is defined for any metric space which have well
defined cut-off $\delta$ covers for all values of $\delta>0$.
Note for example that any simply connected covering space, $X$,
we have $\tilde{X}^{\delta,R}_{cut}=X$ for all $R$ and $\delta$
and so $\tilde{X}^\delta_{cut}=X$ for all $\delta$ and
$CovSpec_{cut}(X)$ is well defined and empty.  The same
thing occurs when $X$ is the standard cylinder, $S^1\times \Bbb{R}$.

By Proposition~\ref{prop4.7}, $CovSpec_{cut}(X)$ is well defined
for all complete locally compact length spaces, $X$, as well.
We expect it is well defined for a much larger class of spaces
but will not be pursuing that investigation ourselves.
See Theorem~\ref{theo4.20} and Examples~\ref{ribbon1} and~\ref{ribbon2}
for other settings where the cut-off covering spectrum is well defined.

Note that by Proposition~\ref{cutcoveringorder}, Theorem~\ref{appendixthm}
and this definition, we have:

\begin{lemma} \label{cutofflowersemiclosed}
The cut-off covering spectrum is a
lower semiclosed set.
\end{lemma}

The following proposition is easy to prove from the definitions.

\begin{prop} \label{boundedX}
If $X$ is a bounded metric space with $D=diam(X)$, then
\be
\tilde{X}^{\delta,R}_{cut} =\tilde{X}^\delta \qquad \forall R\ge D,
\textrm{ and }\tilde{X}^{\delta}_{cut} =\tilde{X}^\delta.
\ee
So $CovSpec_{cut}(X)=CovSpec(X)$.
\end{prop}

Thus the cut-off covering spectrum is really only useful to study
complete length spaces which are not bounded.

In the next subsection we explore the distinction between
these two spectra in general.

\subsection{Relating the various spectra} \label{sect4.3}

The intuitive idea behind the next theorem is that the covering spectrum
can detect any holes that the cut-off covering spectrum sees.

\begin{theo}\label{cutinthecov} \label{thm4.11}
The cut-off covering spectrum of a complete locally compact length space
is a subset of its covering spectrum.
\end{theo}

This follows from
Lemma~\ref{subsets} and Proposition~\ref{bigcupconverse}  which we state and
prove below.

\begin{lem} \label{subsets}
For any basepoint, $x$, in a metric space, $X$,
\be  \label{subset1}
CovSpec^R_{cut}(X,x) \subset CovSpec(X),
\ee
 and
\be  \label{subset2}
CovSpec^{R_1}_{cut}(X,x) \subset CovSpec^{R_2}_{cut}(X,x)
\textrm{ for }R_1<R_2.
\ee
\end{lem}

\Pf  If $\delta \in CovSpec^R_{cut}(X)$, then
$\tilde{X}^{\delta_1, R}_{cut} \neq \tilde{X}^{\delta, R}_{cut}$
for all $\delta_1>\delta$. So there is a nontrivial loop
$\gamma$ which lifts to $\tilde{X}^{\delta, R}_{cut}$
nontrivially and lifts to $\tilde{X}^{\delta_1, R}_{cut}$ trivially.
In particular we can choose
$\gamma$ which lies in a ball of radius $\delta_1$.  Otherwise if all
such loops lift trivially to $\tilde{X}^{\delta, R}_{cut}$
then the covering groups are the same.

If $\delta \notin CovSpec(X)$, then
 $\tilde{X}^{\delta} = \tilde{X}^{\delta_1}$
for some $\delta_1> \delta$. Then $\gamma$ which lifts
trivially to the $\delta_1$ cover, also lifts trivially
to the $\delta$ cover, and must then project trivially
back down to $\tilde{X}^{\delta, R}_{cut}$ nontrivially.
causing a contradiction.

Similarly if $\delta \notin CovSpec^{R_2}_{cut}(X)$, then
 $\tilde{X}_{cut}^{\delta, R_2} = \tilde{X}_{cut}^{\delta_1,R_2}$
for some $\delta_1> \delta$ and we can lift $\gamma$ trivially to both
of these covers which
contradicts that it
lifts to $\tilde{X}^{\delta, R}_{cut}$ nontrivially.
\qed

\begin{prop}\label{bigcup}
If $X$ is a complete locally compact
length space then for any basepoint $x\in X$,
\be  \label{big-cup}
\bigcup_{R>0} CovSpec^{R}_{cut}(X,x)\subset CovSpec_{cut}(X).
\ee
\end{prop}

\Pf If $\delta \in CovSpec^{R_0}_{cut}(X)$, by (\ref{subset2}),
then $\delta \in CovSpec^{R}_{cut}(X)$ for all $R \ge R_0$. So the
covering map
\be
\pi_R:   \tilde{X}^{\delta, R}_{cut} \rightarrow \tilde{X}^{\delta_1,
R}_{cut}
\ee
is nontrivial for all $\delta_1 > \delta$. Then as $R \rightarrow \infty$,
the limit map
\be
\pi: \tilde{X}^{\delta}_{cut} \rightarrow \tilde{X}^{\delta_1}_{cut}
\ee
is nontrivial. So
$\delta \in CovSpec_{cut}(X)$.
Hence $\bigcup_{R>0} CovSpec^{R}_{cut}(X) \subset CovSpec_{cut}(X)$.
\qed

At first one might think that the inclusion in (\ref{big-cup}) is equal.
This is not true.

\begin{example} \label{nobigcup}\label{ex4.2}
Let $X$ be a line with circles attached at the integers $j\neq 0$
of circumference $2\pi r_j$ where $r_j = 1+1/|j|$. Using $0$ as the
base point we have \be CovSpec^R_{cut}(X)= \{\pi+\pi/j: j\in
\Bbb{N}, j+1\le R\} \ee because
 $R$ cut-off $\delta$ covers unravel all loops
such that $j+1 \le R$ and $\pi+\pi/j \ge \delta$.
Taking the Gromov-Hausdorff limit of these
covers we see that
the cut off $\delta$ covers of $X$
unravel all loops $\pi+\pi/j \ge \delta$.
Thus $CovSpec_{cut}(X)$ is the lower semiclosure of
$\{\pi+\pi/j: j\in \Bbb{N}, j+1\le R\}$ which includes the
number $\pi$ because for all $\delta'>\pi$ we have
$\tilde{X}^{\delta'}_{cut}\neq \tilde{X}^{\pi}_{cut}$.
However the union of $CovSpec^R_{cut}(X)$
over all $R>0$ does not include the number $\pi$.
\end{example}

\begin{prop} \label{bigcupconverse}
If $X$ is a complete locally compact
length space then the lower semiclosure of the
union of all $R$ cut-off spectra is the cut-off covering spectrum:
\be
Cl_{lower}\left(\bigcup_{R>0} CovSpec^{R}_{cut}(X)\right)\,\cup\, \{0\}
\,\,=\,\, CovSpec_{cut}(X)\,\cup\, \{0\}.
\ee
\end{prop}

\Pf Take the the lower semiclosure
to both sides of (\ref{big-cup}), since $CovSpec_{cut}(X)$ is lower
semiclosed by Theorem~\ref{appendixthm}, we have
\be
Cl_{lower}\left(\bigcup_{R>0} CovSpec^{R}_{cut}(X)\right) \subset
CovSpec_{cut}(X).
\ee

Now suppose $\delta>0$ is not in the lower semiclosure
of $\bigcup_{R>0} CovSpec^{R}_{cut}(X)$.  Then by Lemma~\ref{appendixopen}
there
exists $\epsilon>0$ such that
\be
[\delta,\delta +\epsilon) \cap \bigcup_{R>0} CovSpec^{R}_{cut}(X)
=\emptyset.
\ee
So for all $R>0$,
\be
[\delta,\delta +\epsilon) \cap  CovSpec^{R}_{cut}(X) =\emptyset
\ee
which implies (by Lemma~\ref{deltaRgap}) that
\be
\tilde{X}^{\delta+\epsilon,R}_{cut}=\tilde{X}^{\delta,R}_{cut}.
\ee
Taking the $R\to\infty$ and the Gromov-Hausdorff limits of
these spaces, we get
\be
\tilde{X}^{\delta+\epsilon}_{cut}=\tilde{X}^{\delta}_{cut}
\ee
which implies that $\delta \notin CovSpec_{cut}(X)$.
\qed

\noindent {\bf Proof of Theorem~\ref{cutinthecov}}: Combining
(\ref{subset1}) with Proposition~\ref{bigcupconverse} the result follows
since $CovSpec (X)$ is a lower semiclosed set.
\qed

\subsection{The length spectrum and the cut-off spectrum}\label{sect4.4}

We now relate the cut-off covering spectrum to the length
spectrum of a complete locally compact length space.
Recall that Example~\ref{nobigcup} is such a space.

\begin{theo} \label{properRcutinL}
If $X$ is complete locally compact length space then
\be
CovSpec^R_{cut}(X) \subset (1/2) L(X).
\ee
That is, if $\delta \in CovSpec^R_{cut}(X)$
then $2\delta \in L(X)$.
\end{theo}

  The assumption that the
space be locally compact is necessary:

\begin{example}\label{boundednotinL}
Let $X$ be the collection of circles of circumeference
$2\pi +2\pi/k$, then
\be
CovSpec_{cut}(X)=CovSpec(X)=\{\pi+\pi/k:\, k\in \Bbb{N}\}\cup \{\pi\}
\ee
while the (1/2) length spectrum of the collection of circles
is all finite sums:
\be
(1/2)Length(X)=\{\sum_{k=1}^\infty a_k \pi(1+1/k): a_k\in \Bbb{N} \}
\ee
which does not include $\pi$.
\end{example}

Before we prove Theorem~\ref{properRcutinL}, we prove the corresponding
proposition which does not require local compactness:

\begin{prop} \label{notproperinL}
Let $X$ be a complete locally compact length space.  If we have
$\delta \in CovSpec^R_{cut}(X,x)$ then there exist
$\delta_j$ decreasing to $\delta$ and
loops, $\sigma_j$, with $L(\sigma_j) < 2\delta_j$,  which are not $\delta$
homotopic to a collection of loops lying outside $\bar{B}(x,R)$.
\end{prop}

\noindent {\bf Proof of Proposition~\ref{notproperinL}:}
Given $\delta \in CovSpec^R_{cut}(X,x)$ we know there
exists $\delta_j$ decreasing to $\delta$ such that
\be
\tilde{X}^{\delta,R}_{cut}\neq \tilde{X}^{\delta_j,R}_{cut}.
\ee
So there exist loops $C_j$ in $X$ which are $\delta_j$
homotopic to loops outside $\bar{B}_x(R)$ but are
not $\delta$ homotopic to such a curve.
Note that $C_j$ is
homotopic to a combination of curves $\alpha \beta \alpha^{-1}$
where $\beta$ lie outside $\bar{B}_x(R)$ or inside $B_p(\delta_j)$.
If all the $\beta$ curves lie outside $\bar{B}_{x}(R)$ then
$C_j$ is $\delta$ homotopic to such curves, so this is impossible.
In fact there must be a  $\beta_j$ which lies in a ball $B_{p_j}(\delta_j)$
which is not $\delta$ homotopic to a collection of loops outside
$\bar{B}_x(R)$.

By Lemma~\ref{balltoL} $\beta_j$ is freely homotopic
to a collection of curves of length $<2\delta_j$.  At
least one of these curves is not $\delta$ homotopic to a
collection of loops outside $\bar{B}_x(R)$
because $\beta_j$ is not.  This is the loop $\sigma_j$.
\qed

We can now add the condition that the space is locally
compact:

\noindent {\bf Proof of Theorem~\ref{properRcutinL}:}
By Proposition~\ref{notproperinL} we have a sequence of curves
$\sigma_j$ in $X$.  Note that $Im(\sigma_j)\cap \bar{B}_x(R)$
is nonempty for all $j$.
Since $L(\sigma_j)<2\delta_j < 4\delta$ for $j$ sufficiently
large
\be
\sigma_j: [0,L(\sigma_j)] \to \bar{B}_x(R+2\delta).
\ee
By the local compactness this closed ball is compact
for $j$ , so we can apply the Arzela-Ascoli theorem to
produce a converging subsequence and a limit loop
$\sigma_\infty$.

It is easy to construct a $\delta$ homotopy from $\sigma_\infty$
to $\sigma_j$ for $j$ sufficiently large so $\sigma_\infty$
is also not $\delta$ homotopic to a loop outside $\bar{B}_x(R)$
and, in particular, not $\delta$ homotopic to a point.
Since
\be
L(\sigma_\infty) \le \liminf_{i\to\infty} L(\sigma_i)
\le \liminf_{i\to\infty} 2\delta_i =2\delta
\ee
we can apply Lemma~\ref{achieve} to say that $\sigma_\infty$
is a closed geodesic and has length $2\delta$ so
$2\delta\in L(X)$.
\qed

Combine this with Proposition~\ref{bigcupconverse}, we get
\begin{coro}\label{coro4.16}
For a complete locally compact length space $X$,
\be
CovSpec_{cut}(X) \subset (1/2) Cl_{lower}(L(X)).
\ee
That is, if $h/2 \in CovSpec^R_{cut}(X)$
then either $h\in L(X)$ or
there exist $h_j \in L(X)$ such that
$h_j$ decrease to $h$.
\end{coro}

Example~\ref{nobigcup} shows that the lower semiclosure is needed here.
In the Section~\ref{sect4.6} we will see that there are complete
length spaces which
are not locally compact, with
well defined
cut-off covering spectra which are not 
in the closure of the length spectrum.

\subsection{Topology and the $CovSpec_{cut}(X)$} \label{sect4.5}

In this section we prove that the cut-off covering
spectrum is empty given certain topological conditions on 
a metric  space $X$:
particularly
Theorem~\ref{thmloopstoinfty} and its converse and
Theorem~\ref{topcross}.  Recall
that the covering spectrum of a simply connected compact
metric space is empty while the cut-off covering spectrum
of a cylinder is empty.
We begin with the loops to infinity property defined in
\cite{So-loops}:

\begin{defn}\label{defloopstoinfty}
Given a metric space, $X$,
a loop $\gamma:S^1 \to X$ is said to have the {\em
loops to infinity} property, if for every compact set $K \subset X$,
there is another loop $\sigma:S^1 \to X\setminus K$
freely homotopic to $\gamma$.

The space $X$ is said to have the loops to
infinity property if all its noncontractible loops
have this property.
\end{defn}

\begin{theo} \label{thmloopstoinfty} \label{thm2.0}\label{theo4.20}
Any metric space $X$ with the loops to infinity property
has a well defined but empty cut-off covering spectrum.
\end{theo}

\Pf
Fix $x_0\in X$
and $\delta>0$.
For every $R>0$ let $K=B_{x_0}(R)$ and for any $g\in \pi_1(X,x_0)$
let $\gamma$ be a representative of $g$ based at $x_0$.  So there
exists $\beta$ freely homotopic to $\gamma$ outside $K$ which means
there is a curve $\alpha \circ \beta\circ \alpha^{-1}$ which
represents $g$ such that $\beta \subset X\setminus \bar{B}_p(R)$.
So every $g \in \pi_1(X,x_0)$ is in the covering group of
$\tilde{X}^{\delta,R}_{cut}$, which means
$\tilde{X}^{\delta,R}_{cut}=X$.  Taking the limit $R\to\infty$
we get $\tilde{X}^{\delta}_{cut}=X$ for all $\delta$ so the
cut-off covering spectrum is well defined but trivial.
\qed

This theorem is applied to complete manifolds with
nonnegative Ricci curvature in Theorem~\ref{Riccicutoff}.
Such manifolds have only one end.

Recall that a length space $X$ is said to have $k$ ends if for all
sufficiently large compact sets $K$, $X\setminus K$ has $k$ path
connected components.

A length space is semilocally simply connected if every point
has a neighborhood around it such that any curve in that neighborhood
is contractible.  A Riemannian manifold is semilocally simply connected.

\begin{theo} \label{thm-one-end}\label{thmoneend}
Let $X$ be a complete, locally compact and semilocally simply connected
length space
with an empty cut-off covering spectrum, then any curve in $X$
is homotopic to a product of
curves which have the loops to infinity property.
If in addition $X$ has only one end
then $X$ has the loops to infinity property.
\end{theo}

Example~\ref{extwoend}, right below, demonstrates the necessity
of the one end hypothesis while
Example~\ref{ribbon2} demonstrates the necessity of the local
compactness condition.

\Pf If the cut-off covering spectrum is empty then
 $\tilde{X}^\delta_{cut}=X$ for all $\delta >0$ and, by
Proposition~\ref{cutcoverexists},
$\tilde{X}^{\delta,R}_{cut}$ is between these two spaces, so
it is isometric to $X$ as well.  Thus for all $\delta>0$ and
for all $R>0$, the fundamental group of $X$ is generated by elements
of the form $\alpha\circ\beta\circ\alpha$ where $\beta$ is either
in a ball of radius $\delta$ or in $X\setminus \bar{B}_{x_0}(R)$.

Choose any nontrivial loop $\gamma$
and any compact set $K\subset X$.  Take $R>0$ large enough that
\be
K\cup Im(\gamma)\subset B_{x_0}(R/2).
\ee
Since $X$ is  complete, locally compact and semilocally simply connected, we
can
take $\delta>0$ small enough that balls of radius $\delta$
in $B_{x_0}(R)$ are semilocally simply connected so that any loop
$\beta$ in such a ball is contractible.  Thus
$[\gamma]\in \pi_1(X,x_0)$ must be generated by loops
of the form $\alpha\circ\beta\circ\alpha^{-1}$ where
\be
\beta \in X\setminus \bar{B}_{x_0}(R)\subset X \setminus K.
\ee
When $X$ has only one end, the set $X\setminus K$ is path
connected, thus the various $\beta$ used to generate $X$
can be connected via new paths $\alpha \in X \setminus K$
to a point $x_1 \in X\setminus K$.  Thus we have constructed
$\sigma \in X \setminus K$ which is freely homotopic to
$\gamma$.
\qed

\begin{example}\label{extwoend}
One end is necessary as can be seen by taking the length
space $X$ formed by joining two closed half cylinders
at a point.  The loop $\gamma$ running around a figure eight
which goes once around each cylinder, does not
have the loops to infinity property.  It is generated by
2 different loops $\beta_j$ each of which goes to infinity
in a different direction.
 This can be made smooth
by taking the connected sum of two manifolds that are not
simply connected that have only one end each, like
Nabonnand's example \cite{Nab}.
\end{example}

\begin{theo}\label{topcross}
If a complete length space $X$ is homeomorphic to
the product of complete locally compact length spaces, $M\times N$, then
$X$ has the loops to infinity property and
$CovSpec_{cut}(X)=\emptyset$ if
either of the following holds:

i) both $M$ and $N$ are noncompact

ii) $M$ is noncompact  and $CovSpec(M)=\emptyset$.
\end{theo}

\Pf
Let $C$ be a loop in $X$, so $C=(a,b)$ where
$a$ and $b$ are closed loops in $M$ and $N$
respectively.  $C$ is freely homotopic to
$(a,b_0)$ followed by $(a_0,b)$ where
$b_0=b(0)$  and $a_0=a(0)$.

In both cases $M$ is complete and noncompact, so there
exists $p_j\in M$ which diverge to infinity and
there exist
minimal paths $\sigma_j$ from any fixed point $p_0$
to $p_j$.  If $b$ is a loop in $N$, then
$(p_0,b)$ is freely homotopic to $(p_i,b)$ via
$(\sigma_j,b)$.  Any compact $K\subset X$,
is a subset of the image of $K_M\times K_N$
where $K_M$ is compact in $M$, taking $p_j\in M\setminus K_M$
we have $(p_j,b)$ outside $K$.  Thus $(p_0,b)$
has the loops to infinity property.

In case i, $N$ is also noncompact so both $(a,b_0)$
and $(a_0,b)$ have the loops to infinity property.
So any loop $C$ in $X$ is a combination of curves
with the loops to infinity property and we just
apply Theorem~\ref{thmloopstoinfty}.

Before we begin case ii we note that:
{\em if $a$ has the loops to infinity property
then so does $(a,b_0)$}.  This is seen by
taking the homotopies $h_i$ from $a$ to $a_i$
that diverge to infinity.   Mapping them
to $X$, we get homotopies
$(h_i,b_0)$ from $(a,b_0)$ to $(a_i,b_0)$.
So for any compact set $K\subset X$,
we have $K\subset K_M\times K_N$ where
$K_M$ is compact.  So
we can choose $a_i$ in $M \setminus K_M$
and have $(a_i,b_0)$ outside $K$.

In case ii, we don't have ray in $N$ for the loops in $M$, but
$CovSpec_{cut}(M)=\emptyset$.  Applying
Theorem~\ref{thm-one-end}, we see
that the loop $a$ in $M$ is freely homotopic
to a combination of loops which have the loops
to infinity property.  Thus $(a,b)$ is freely
homotopic to a combination of loops $(a_i,b)$
each of which is homotopic to $(a_i(0),b)$
following $(a_i, b(0))$.  Each $(a_i,b(0))$ has
the loops to infinity property via the loops to
infinity property of each $a_i$.  As in case i),
each $(a_i(0), b)$ has the loops to infinity property
via rays in $M$ based at $a_i(0)$.  So
$CovSpec_{cut}(X)=\emptyset$ here as well.
\qed

\begin{coro}\label{corotop}
If $X$ is a complete noncompact locally compact length space
homeomorphic to $M\times \Bbb{R}$
then 
\be
CovSpec_{cut}(X)=\emptyset.
\ee
\end{coro}

\subsection{Pulled Ribbon Spaces} \label{sectribbon}

In this section we construct examples of metric spaces
with well defined cut off covering spectra demonstrating that
Theorem~\ref{properRcutinL} does not hold without the assumption of
local compactness.  We call the method of construction the
``pulled ribbon construction''.  It is similar to an idea
in Burago-Burago-Ivanov called a ``pulled string'', where
a collection of points lying on a path in a space is
identified creating a new complete length space.  Their construction
is called the pulled string construction because it looks
something like a cloth which has had a thread pulled tight.
In our case we first attach a ribbon along
the line in the space and then we pull a string on the
opposite edge of the ribbon.

We will make our construction precise.  Those who wish
to understand their construction may consult \cite{BBI}.

\begin{defn}
The pulled ribbon space is a space $Y= \Bbb{R}\times [0,1]$
with the lower boundary $\Bbb{R}\times \{0\}$ identified with
a point and endowed with the induced length structure.  This is
the same as saying that the metric on $Y$ is
\be
d_Y( (r_1,s_1), (r_2,s_2))=
\min \{ \sqrt{ (r_1-r_2)^2+(s_1-s_2)^2 }, s_1+s_2 \}
\ee
This is a quasi metric and becomes a metric when we make
the identification $(r_1,0)=(r_2,0)$.
\end{defn}

Note that $Y$ is the suspension of a line.  There is
a geodesic $\gamma:\Bbb{R} \to Y$ which runs along the
``top edge'': $\gamma(r)=(r,1)$.  This geodesic is
not a line in the induced length structure.

\begin{prop}
The pulled ribbon space is a bounded
complete length space which is not locally compact.
\end{prop}

\Pf
It is bounded because $Y\subset \bar{B}_{y}(1)$ where
$y$ is the special identified point.  It is not compact because
the sequence of points $(2j,1)$ are all a distance $2$ apart
from each other.  It is a complete length space because between
any pair of points we can find a minimal geodesic between them:
it is either the line segment in the strip or a pair of vertical
lines dropping from the points to the common point.  Given any
Cauchy sequence $(r_i,s_i)$ in the induced length structure
either the sequence converges in the standard metric on the strip
or it approached the bottom edge which is the common point.
\qed

\begin{defn}
Given a manifold $M$ with a line $\gamma:\Bbb{R}\to M$,
we say that we attach a pulled ribbon to $M$ creating
a space, $M_\gamma$, if we attach the pulled ribbon
so that its top edge is identified with the line.  Then
we endow $M_\gamma$ with the induced length metric.
\end{defn}

Note that in this induced length metric the original line
$\gamma$ is no longer a line and is now bounded.  However,
unlike spaces with a pulled thread, a space with an attached
pulled ribbon keeps its topology.  In fact:

\begin{prop} \label{lessthan1}
If $x,y\in M$ and $d_M(x,y)<2$ then after
adding the pulled ribbon to $M$, we do not change the
distance between $x$ and $y$.
\end{prop}

\Pf
If the distance between $x$ and $y$ has been
shortened then there is a path from $x$ to $y$
of length $<2$ which passes into the ribbon.
However, such a short path could not reach the
far edge of the ribbon, and so it's length
is determined by the Euclidean structure on the
ribbon and it would be shorter if it did not enter
the ribbon at all.
\qed

\begin{coro}
Given a space $M$ and a map $f: X \to M$  then
$f$ is continuous from $X$ to $M$ iff
$f$ is continuous from $X$ to $M_\gamma$ with
the ribbon attached.
\end{coro}

\begin{coro} \label{ribbonsimply}
$M$ is simply connected iff $M_\gamma$
is simply connected.  Furthermore,
$M$ is semilocally simply connected iff
$M_\gamma$ is semilocally simply connected.
\end{coro}

The following example demonstrates the necessity of
the local compactness condition in Theorem~\ref{thmoneend}.

\begin{example}\label{ribbon1}
Let $M^2$ be the cusped manifold,
\be
\mathbb R\times_f S^1,
\ee
where $f(r)=2Arctan(-r)+\pi$ and $\gamma$ be
any line in this space.

Then $M_\gamma$ is not simply connected but is semilocally simply connected
by Corollary~\ref{ribbonsimply}.
It is a bounded space which is not locally compact.
Note that $\tilde{M}^{\delta,R}_{cut}=M$ for every value of $R,\delta>0$.
So $CovSpec_{cut}^R(M))$ is trivial and $CovSpec_{cut}(M)$ is
well defined and trivial. 

Adding a halfline attached at any point, would just create a space which
 does not have the loops to infinity property and has one end.
Nevertheless,
\be
CovSpec_{cut}(M_\gamma)=CovSpec_{cut}^R(M_\gamma)=CovSpec(M_\gamma)=\emptyset.
\ee
\end{example}

Note that we just pulled a thread in this example identifying
a line to be a point and using the induced length structure, the space
would become simply connected: loops shrinking along the cusp
would in fact converge to the identified point which is the line.
The loops in this example do not have a converging subsequence because
they are always a fixed distance away from the line.

It should be noted that the double suspension of the Hawaii Ring
is a compact space which is not simply connected and yet it
is its own universal cover so its covering spectrum is empty as well
\cite{Sp}.

We next demonstrate that  local compactness is necessary in
Theorem~\ref{properRcutinL}.

\begin{prop}\label{ribbonlength}
If the length spectrum of $M$ is empty and there
are no geodesics starting and ending perpendicular
to the line $\gamma$, then the length
spectrum of $M_\gamma$ is also empty.
\end{prop}

\Pf
Suppose on the contrary that there is a closed geodesic
$\sigma:S^1 \to M_\gamma$.  If its image lies in $M\subset M_\gamma$,
then it is also a closed geodesic in $M$ by Proposition~\ref{lessthan1}.
Since $M$ has an empty length spectrum this cannot be the case, so
its image must intersect with the ribbon.  The image of a closed
geodesic cannot lie completely within the ribbon, because there are
no closed geodesics formed using Euclidean line segments.  So the
geodesic $\sigma$ must enter and leave the ribbon.  The only way
$\sigma$ could turn around is if it passes through the far edge
and comes back.  Thus the geodesic must be vertical and must
intersect the line $\gamma$ vertically.  So the part of $\sigma$
which lies in $M$ contradicts the hypothesis.
\qed

This next example demonstrates the necessity of the local compactness
condition in Theorem~\ref{thmoneend}:

\begin{example}\label{ribbon2}\label{nolift}
Let $M$ be asymptotically cylindrical
\be
\mathbb R\times_f S^1
\ee
where $f(r)=2Arctan(-r)+2\pi$ and $\gamma$ be
any line in $M$, then $M_\gamma$
is a bounded complete geodesic length space
with diameter $D\le 2+2\pi$.  It has 
an empty length spectrum by
Proposition~\ref{ribbonlength}.

Note that for $\delta \le \pi$, $\tilde{M}^\delta=\tilde{M}$
which is a complete length space that does not have minimizing
geodesics joining every pair of points (particularly the lifts
of the pulled point).  Note also that for $R\ge D$ and 
$\delta\le \pi$, $\tilde{M}^{R,\delta}_{cut}=\tilde{M}$
while for $R \ge D$ and $\delta>\pi$, $\tilde{M}^{R,\delta}_{cut}=M$.
Thus the cut-off $\delta$ covers are defined and so is
the cut-off covering spectrum.  Furthermore
\be
CovSpec_{cut}(M_\gamma)=CovSpec_{cut}^D(M_\gamma)=CovSpec(M_\gamma)=\{\pi\}.
\ee
\end{example}

\subsection{Localizing the $R$ cut-off covering spectrum} \label{sect4.6}

In this section we show that one can compute
$CovSpec_{cut}^R(X,x)\cap[0,D]$ of a metric space $X$ 
using only the information contained in $B(x,r)$
when $r$ is taken sufficiently large [Prop~\ref{isomballs}].
In fact we give a precise estimate on $r$ independant
of $X$ which will allow us to stufy sequences of
spaces.

Note that there is a complete hyperbolic manifold $M$
of constant sectional curvature $-1$ such that for any
$r$, there exists a contractible curve lying in $B(p,1)$ which is not
homotopically trivial in $B(p,r)$ \cite{BoMe}\cite{Po}.
In other words, the homotopies required to contract these loops
to a point extend further and further out in $M$.
A simpler example with this property is formed by
taking the Hawaii Ring with circles of circumference $1/k$
and attaching a cylinders of length $k$ to the $k^{th}$
circle and then capping off the cylinder.  This is a simply
connected space none of whose balls about the basepoint
are simply connected.

\begin{figure}[htbp]
\includegraphics[width=4.5in ]{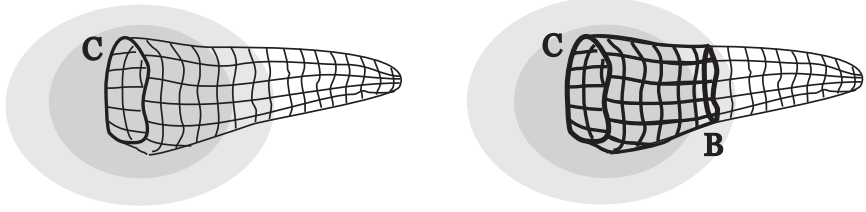}
\caption{Here $C$ is $\delta$ homotopic to
a single $\beta=B$ outside $\bar{B}_x(R)$.}
\label{cutnet}
 \end{figure}

The covering spectrum of these spaces could not be computed
using a localization process like the one we obtain here
for the $R$ cut-off covering spectrum.  It is crucial that
we can chop  off homotopies as in Figure~\ref{cutnet}
when computing the $R$ cut-off covering spectrum.

Recall the definition of $\delta$ homotopy in Definition~\ref{defdelhom}
and Lemma~\ref{equivdefdelhom} and Defn~\ref{extdefdelhom}.
Now we define:

\begin{defn}
Two loops $\gamma_1, \gamma_2$ in a metric space, $X$,
are $R$-cutoff $\delta$ homotopic in $X$ if
$\pi_{\delta, R} (\gamma_1) = \pi_{\delta,R} (\gamma_2)$,
where $\pi_{\delta,R}: \pi(X) \rightarrow \pi(X)/\pi(X, \delta,R)$.
\end{defn}

It is not hard to see from the definition of the $R$ cut-off
$\delta$ cover that we have the following simpler description
which will allow us to apply the lemmas from the section on
$\delta$ homotopies to study this new kind of homotopy:

\begin{lemma} \label{Rcutrelate}
Let $A$ be a subset of a metric space $X$.
A loop $\gamma$ is $R$-cutoff $\delta$ homotopic to a point
in $A$ iff it is $\delta$ homotopic in $A$ to a collection of
loops $\beta_j$ lying outside $\bar{B}_p(R)$.
\end{lemma}

Our next lemma will be useful for localizing the $\delta$
homotopies so that we can use compactness to control them.

\begin{lemma}\label{cuthomotopy}
Given $\delta >0, R> 0$, and a loop $C$ in $B(x,R+2\delta)$
in a metric space, $X$,
if $C$ is $\delta$ homotopic in $X$
to a collection of curves $\alpha\beta\alpha^{-1}$
where $\beta$ are in $\delta$-balls or outside $\bar{B}(x,R)$, then
$C$ is $\delta$-homotopic in $B(x,R+2\delta)$
to a collection of curves $\alpha\beta\alpha^{-1}$
where $\beta$ are in $\delta$-balls or outside $\bar{B}(x,R)$.  So $C$ is
$R$-cutoff $\delta$ homotopic to a point in $\bar{B}(x,R+2\delta)$.
\end{lemma}

See Figure~\ref{cutnet} where the darker balls are $B(x,R)$
and the lighter balls are $B(x,R+2\delta)$.

\Pf
This proof follows from Lemma~\ref{squaregridlem}
where our set $A=X$ and $B=B(x, R+2\delta)$
so we see that $C$
is $\delta$ homotopic in $B$ to a collection
of curves $\gamma_1, \gamma_2,...\gamma_j$ such that
each $\gamma_j$ lies in
$B$ and the tubular neighborhood
\be
T_{2\delta}(A\setminus B).
\ee
In particular the $\gamma_i$ lie outside $\bar{B}(x,R)$.
Thus by Lemma~\ref{Rcutrelate},
$C$ is a curve which is $R$-cutoff $\delta$-homotopic in $B(x,R+2\delta)$
to a collection of curves $\alpha\beta\alpha^{-1}$
where $\beta$ are in $\delta$-balls or outside $\bar{B}(x,R)$.
\qed

Using Lemma~\ref{cuthomotopy} we have the following relation
between the $R$-cutoff spectrums of balls and the total space
which will be very useful later.

\begin{lemma}\label{uniflocaltoglobal}
If $X$ is a metric space, then
\[ \delta \in CovSpec^{R}_{cut}(B(x,r))
\]
for some $r > 3(R+2\delta)$ then
\[
\delta \in CovSpec^{R}_{cut}(X).
\]
When $X$ is a length space we can use the induced
length metric on $B(x,r)$ rather than the restricted
metric and still have the same result.
\end{lemma}

\Pf
We prove the contrapositive.  Assume $\delta \notin CovSpec_{cut}^R(X)$.
By the definition, there exists $\delta'>\delta$
such that
\be
\tilde{X}_{cut}^{\delta,R} = \tilde{X}_{cut}^{R,\delta'}.
\ee
This means that any curve $C$
whose image lies in a ball of radius $\delta'$  is $\delta$-homotopic in $X$
to a
path created as a combination of $\alpha \beta \alpha^{-1}$
where $\beta$ are either in a ball of radius $\delta$ or
lie outside $\bar{B}_x(R)$.

If the image of $C$ lies in $B(x, R+2\delta')$, then by
Lemma~\ref{cuthomotopy}
$C$ is a curve which is $\delta$-homotopic in $B(x,R+2\delta')$
to a collection of curves $\alpha\beta\alpha^{-1}$
where $\beta$ are in $\delta$-balls or outside $\bar{B}(x,R)$.
When $r> 3(R+2\delta)$ and $\delta'$ is close to $\delta$, the metric on
$B_x(R+2\delta')$ restricted from the induced length metric on $B_x(r)$
agrees with its metric
restricted from $X$, we have that any such curve $C$
lifts closed to $\tilde{B}(x,r)_{cut}^{\delta,R}$.

Now any curve $\sigma$ which lifts closed to
$\tilde{B}(x,r)_{cut}^{R,\delta'}$ is homotopic in $B(x,r)$ to
a collection of $\alpha  \beta \alpha^{-1}$
where $\beta$ are now either in a ball of radius $\delta'$
or outside $\bar{B}_x(R)$.  Note that any $\beta$ which
lies outside $\bar{B}_x(R)$ lifts as a closed loop to
$\tilde{B}(x,r)_{cut}^{\delta,R}$.  Those $\beta$ which
pass within $\bar{B}_x(R)$ and fit in a ball of radius
$\delta'$, the images must be contained in  $B(x, R+2\delta')$, therefore
also
lifts closed to $\tilde{B}(x,r)_{cut}^{\delta,R}$.
Since $\sigma$ is homotopic in $B(x,r)$ to a combination
of curves which lift as closed loops to
$\tilde{B}(x,r)_{cut}^{\delta,R}$, then $\sigma$ must do the same.
Thus by the Curve Lifting property (c.f. \cite{Massey} page 123)
and Lemma~\ref{cutcoveringorder} we see that
\be
\tilde{B}(x,r)_{cut}^{\delta',R}=\tilde{B}(x,r)_{cut}^{\delta,R}
\ee
and so $\delta \notin CovSpec_{cut}^R(B(x,r))$.
\qed

In the opposition direction we have
\begin{lemma}\label{globaltolocal}
If $X$ is a metric space with
\[
\delta \in CovSpec^{R}_{cut}(X)
\]
then for all $r\ge 3(R+2\delta)$,
\[
\delta \in CovSpec^{R}_{cut}(B(x,r)).
\]
Again when $X$ is a length space we may either
give $B(x,r)$ the induced length metric or
the restricted metric.
\end{lemma}

\Pf If
\be
\delta \in CovSpec^{R}_{cut}(X),
\ee
then, for  $\delta_i$ decreasing to $\delta$,  we have loops $C_i$ lies
inside $\delta_i$ balls of $X$ which can not be represented in $X$ by loops
lying inside $\delta$-balls of $X$  or loops in $X\setminus \bar{B}(x,R)$.
Since $C_i$ lies inside $\delta_i$ balls and is not in $X\setminus
\bar{B}(x,R)$ it must be in $B(x,R+2\delta_i)$. Since $r \ge 3(R+2\delta)$
the balls in $B(x,R+2\delta_i)$ are same for $B(x,r)$ and $X$. Therefore
$C_i$ lies inside $\delta_i$ balls of $B(x,r)$. Since $C_i$ can not be
represented in $X$ by loops lying inside $\delta$-balls of $X$  or loops in
$X\setminus \bar{B}(x,R)$  it can not be represented by loops lying inside
$\delta$-balls of $B(x,r)$ or loops in $B(x,r) \setminus \bar{B}(x,R)$   in
$B(x,r)$. This shows that $\delta \in CovSpec^{R}_{cut}(B(x,r))$.
\qed

An immediate consequence of these two lemmas is the following:

\begin{prop}\label{isomballs}
Given two metric spaces $X$ and $Y$
with isometric balls, $B(x,r)=B(y,r)$, then
\be
CovSpec_{cut}^R(X,x)\cap [0,D] =CovSpec_{cut}^R(Y,y)\cap [0,D]
\ee
whenever $3(R+2D) \le r$.  These balls may have restricted metrics
or induced length metrics when $X$ and $Y$ are length spaces.
\end{prop}

\Pf
If $\delta \in CovSpec_{cut}^R(X,x)\cap [0,D]$,
then $\delta\le D$ so apply Lemma~\ref{globaltolocal}
and have
\be
\delta \in CovSpec^{R}_{cut}(B(x,r))=CovSpec^{R}_{cut}(B(y,r)).
\ee
Then apply Lemma~\ref{uniflocaltoglobal} gives the result.
\qed

Remember that the $R$ cut-off covering spectrum of a capped
cylinder and a cylinder are both empty regardless of basepoint
while the ordinary covering spectrum of the cylinder is nonempty.

Without restricting to a uniform $[0,D]$,
the $R$ cut-off covering spectrum will not match.  This can
be seen in the following example:

\begin{example}
Let $X_s$ be a unit interval with $x_s$ on one end and
a circle  of circumference $2\pi s$ on the other end.
Let $Y$ be a unit interval with $y$ at one end and
two half lines at the far end.  Taking $R=2$ and
$s>1$ we have $CovSpec_{cut}^R(X_s,x_s)=\{\pi s\}$ and
$CovSpec_{cut}^R(Y,y)=\emptyset$.  Yet for any $r$
we have $B(x_s,r)$ isometric to $B(y,r)$ for $s>r$.
\end{example}

\subsection{Varying $R$ in the $R$ cut-off covering spectra} \label{sect4.7}

In the next section on the Gromov-Hausdorff convergence
of metric spaces and the cut-off covering spectra we
need to relate the $R$ cut-off covering spectra for various
values of $R$.

\begin{prop} \label{Rflex}\label{prop4.26}
Given a metric space $X$ and $R_0<R_1$,
\[
\delta\in CovSpec^{R_1}_{cut}(X) \setminus
CovSpec^{R_0}_{cut}(X)
\]
implies
\[ \label{changeR}
\tilde{X}^{\delta, R_1}_{cut} \to  \tilde{X}^{\delta, R_0}_{cut}.
\]
is nontrivial.
\end{prop}

\Pf  If $\delta \in CovSpec^{R_1}_{cut}(X)$, then
$\tilde{X}^{\delta_i, R_1}_{cut} \neq \tilde{X}^{\delta, R_1}_{cut}$
for all $\delta_i>\delta$. So there is a nontrivial loop
$\gamma_i$ which lifts to $\tilde{X}^{\delta, R_1}_{cut}$
nontrivially and lifts to $\tilde{X}^{\delta_i, R_1}_{cut}$ trivially.
Since $R_1$ is the same for both covering spaces,
we can choose
$\gamma_i$ which lies in a balls of radius $\delta_i$.  Otherwise if all
such loops lift trivially to $\tilde{X}^{\delta, R}_{cut}$
then the covering groups are the same.

Suppose $\delta \notin CovSpec^{R_0}_{cut}(X)$,
then for $i$ sufficiently large,
\be
\tilde{X}^{\delta_i, R_0}_{cut} = \tilde{X}^{\delta, R_0}_{cut}.
\ee
Since $\gamma_i$ lies in a ball of radius $\delta_i$
it lifts trivially to the first cover, and thus also the second.
So we have a nontrivial covering:
\be
\tilde{X}^{\delta, R_1}_{cut} \to  \tilde{X}^{\delta, R_0}_{cut}.
\ee
\qed

In the next proposition we assume our space $Y$ is compact.  To apply
this proposition to complete noncompact spaces $X$ which are only locally
compact
we will use our localization results from the last section.

\begin{prop} \label{iflex} \label{prop4.27}
If $Y$ is a compact  length space,
and  $CovSpec_{cut}^{R_i}(Y) \cap [\delta_1,\delta_2)=\emptyset$
for a sequence of $R_i$ decreasing to $R_1$, then
for $R_i$ sufficiently close to $R_1$ we have
\be \label{concit}
\tilde{Y}_{cut}^{\delta_1,R_i} \to \tilde{Y}_{cut}^{\delta_2,R_1}
\ee
is trivial. In particular, without any assumption on
the spectrum, we have
\be \label{concit2}
\tilde{Y}_{cut}^{\delta,R_i} \to \tilde{Y}_{cut}^{\delta,R_1}
\ee
is trivial whenever $R_i$ is sufficiently close to $R_1$.
\end{prop}

Combining this with Proposition~\ref{subsets}, we only
need to assume there exists $R_2>R_1$ with
$CovSpec_{cut}^{R_2}(Y) \cap [\delta_1,\delta_2)=\emptyset$
to conclude (\ref{concit}).
In fact, by Proposition \ref{bigcup} we could assume
$CovSpec_{cut}(Y) \cap [\delta_1,\delta_2)=\emptyset$
and draw the same conclusion.
See Theorem~\ref{covcutconv} for an application of this
proposition.

\Pf Assume on the contrary that
$\tilde{Y}^{\delta_1,R_2} \to \tilde{Y}^{\delta_2,R_1}$
is not trivial for all $R_2 > R_1$.  So there is a $\gamma$ which lifts
trivially to the latter cover, but not to the
first.  In particular we can either choose $\gamma$
to lie inside a ball of radius $\delta_2$, or
outside $\bar{B}_p(R_1)$.

In the first case,
$\gamma$ lifts trivially to $\tilde{Y}^{\delta_2,R_2}$
which implies $CovSpec_{cut}^{R_2}(Y) \cap [\delta_1,\delta_2)$
is nonempty.

In the second case $\gamma$ lies outside
 $\bar{B}_p(R_1)$ and is not $\delta_1$
homotopic to a loop outside $\bar{B}_p(R_2)$.
In particular $l(\gamma)\ge 2\delta_1$
and $\gamma$ is not $\delta$ homotopic
to a loop outside $\bar{B}_p(R_2)$ for
any $\delta\le \delta_1$.

Suppose we take $R_2=R_i$ decreasing to $R_1$
and have nontrivial covers.  So we get a sequence
of $\gamma_i$, each $\gamma_i$ lies outside
 $\bar{B}_p(R_1)$ and is not $\delta_1$
homotopic to a loop outside $\bar{B}_p(R_i)$.

Note that $Y \setminus \bar{B}_p(R_1)$ is precompact.  It is still a
precompact length space if we give it the induced length structure
(c.f. \cite{BBI}).  So there exists some finite number $N$ such that
it can be covered by at most $N$ balls of radius $\delta_1/5$. Note
that balls in the induced length metric are smaller than those in
the metric on $Z$, so $\gamma_i$ is also not $\delta_1$ homotopic in
the space $Z=Y \setminus \bar{B}_p(R_1)$ to a loop outside
$\bar{B}_p(R_i)$. For the rest of the proof we will use the induced
length metric on $Z$ when referring to the $\delta_1$ homotopies..

Applying Lemma~\ref{lengthbound}, we see that
we can always find a $\gamma_i$ in $Z$
with $L(\gamma_i)\le N\delta_1$
which is not $\delta_1$ homotopic to
a loop outside $\bar{B}_p(R_i)$.

Since the $\gamma_i$ have length bounded above uniformly and
since $Y$ is compact, by Arzela Ascoli we have a subsequence
which converges to some $\gamma_\infty$.  Note that $\gamma_\infty$
need not be located outside $\bar{B}_p(R_1)$, so instead of
relating $\gamma_i$ to $\gamma_\infty$, we will use the fact that
$\gamma_i$ must be a Cauchy
sequence in $Z$.  That is , there exists $N'$ sufficiently
large such that $\gamma_i$ are $\delta_1/2$ homotopic to $\gamma_j$
for all $i,j \ge N'$. Fix this $N$, note that
$\gamma_N$ lies outside the closed ball
$\bar{B}_p(R_1)$ and $R_j$ are decreasing to
$R_1$, so $\gamma_N$ is outside
$\bar{B}_p(R_j)$ for $j$ sufficiently
large.  This contradicts $\gamma_j$ is not $\delta_1$
homotopic to a loop outside $\bar{B}_p(R_j)$.
\qed

Note that the compactness here is essential as the following example
shows.

\begin{example}
Let $Y$ be the Hawaii ring with circles of circumference $2\pi \pm \frac
{\pi}{j}$, $\gamma_j$, all attached at a point. Take $\delta =\pi/2, \ R_i =
(1+1/i)\pi, \ R_1 =\pi$, then the cover
$\tilde{Y}_{cut}^{\delta,R_i} \to \tilde{Y}_{cut}^{\delta,R_1}$ is
nontrivial for all $i$. This $Y$ is not a compact length space.
\end{example}


\sect{Gromov-Hausdorff Convergence} \label{sect5}\label{Sect5}

In \cite{SoWei3}, we proved that when compact spaces $M_j$
converge to a compact limit $M$ in the Gromov-Hausdorff sense
then $CovSpec(M_i)\cup \{0\}$ converges to $CovSpec(M) \cup \{0\}$
in the Hausdorff sense as subsets of the real line.
In particular, if $M_j$ are simply connected, then the limit space has
an empty covering spectrum and is its own universal cover.

In the next subsection we provide examples demonstrating that
we do not get such a strong result when the spaces are noncompact.
In fact the limit space of simply connected $M_i$
might be a cylinder [Example~\ref{cylinder}].

In the subsequent sections we prove the continuity of the
cut-off covering spectra [Theorem~\ref{covcutconv}].
In particular the limit of simply connected manifolds will be
seen to have an empty cut-off covering spectrum [Corollary~\ref{limsimp}].

First recall the definition of Gromov Hausdorff distance:

\begin{defn}
Given compact length spaces $X_i$ and $Y$ we say
$X_i$ converges to $Y$ in the Gromov Hausdorff sense if there
exists $\delta_i$ Hausdorff approximations
$f_i: X_i \to Y$ such that
\be
|d_Y(f_i(x_1),f_i(x_2))-d_{X_i}(x_1,x_2)|<\delta_i
\ee
and $Y \subset T_{\delta_i}(f(X_i))$ with $\delta_i \to 0$.
Note that once this is true there are also $\delta'_i$
Hausdorff approximations from $Y$ to $X_i$ with $\delta'_i\to 0$.
\end{defn}

When complete noncompact spaces are said to
converge in the Gromov-Hausdorff sense, they are considered as pointed
spaces.  We write $(X_i,x_i)$ converges in the pointed Gromov-Hausdorff
sense to $(X,x)$ when for every $R>0$, the closed balls with
the restricted metric $\bar{B}_{x_i}(R)\subset X_i$ converge to balls in
the limit space $\bar{B}_{x}(R) \subset X$.

Gromov's Compactness
Theorem says that a sequence of complete locally compact metric spaces, $X_k$,
converges in the pointed Gromov-Hausdorff sense iff they
 are uniformly locally compact in the sense that the number of
disjoint balls of radius $\epsilon$ lying in a ball of radius $R$
is uniformly bounded $N(\epsilon, R, X_k) \le N(\epsilon,R)$. 
Crucial here is that the balls of increasing radius $R$ do piece together
to form a complete locally compact limit space which is also a length space
when the $X_k$ are length spaces.  \cite{Gr}
However, one must keep in mind that
the balls can converge at different rates.  The next section depicts
a few examples where this aspect of the pointed Gromov-Hausdorff
convergence is crucial.

\subsection{Examples}\label{sect5.1}

First recall that even for the covering spectrum on
compact spaces it is possible for a sequence of spaces to
become simply connected in the limit:

\begin{example}\label{shrinkto0} \label{ex5.1}
Let $M$ be a simply connected surface and let $X_k$ be created by
adding a small handle onto $M$, such that the handle fits inside a
ball of radius $1/k$. These $X_k$ converge to $M$ as $k \to \infty$.
Note that both $CovSpec(X_k)=\{\delta_k\}$ and
$CovSpec_{cut}(X_k)=\{\delta_k\}$ while
$CovSpec(M)=\emptyset$.
\end{example}

Nevertheless in \cite{SoWei3} we proved that the difficulty seen
here was the only cause for a lack of continuity in the covering
spectrum.  We proved for compact $M_i$ converging to compact limits
$Y$,  then if $\lambda_j \in CovSpec(M_i)$ converge to $\lambda>0$
then $\lambda \in CovSpec(Y)$ and if $\lambda\in CovSpec(Y)$ there
exists $\lambda_j\in CovSpec(M_i)$ such that $\lambda_j \to
\lambda$.  In particular, if the $M_i$ are simply connected, then $Y$ has an
empty covering spectrum.

Without the assumption of compactness, however,
we can have simply connected manifolds which have
a limit with a nonempty covering spectrum:

\begin{example}\label{cylinder} \label{ex5.2}
Let $M$ be a capped off cylinder and let $p_i\in M$ diverge
to infinity.  Then the sequence $(M_i,p_i)$ converges in the
pointed Gromov-Hausdorff sense to a cylinder because the cap
has disappeared off to infinity.

Thus we have a sequence $M_i \to Y$ such that $CovSpec(M_i)=\emptyset$
but $CovSpec(Y)=\{\pi\}$.  Now by Propositions~\ref{subsets}
and~\ref{bigcup}
$CovSpec_{cut}(M_i)=CovSpec_{cut}^R(M_i)=\emptyset$ as well.
Since a cylinder has the loops to infinity property
$CovSpec_{cut}(Y)=CovSpec_{cut}^R(Y)$ are also empty.
\end{example}

While in the above example, the limit gained an element in
its covering spectrum due to longer and longer homotopies, it is also
possible to gain an element in the covering spectrum without
changing the topology of the space:

\begin{example} \label{widening}\label{ex5.3}
{\em We construct an example where an element of the covering spectrum
appears in the limit.}
As in Example~\ref{slipex},
let
$M^2$ be the warped product manifold $\Bbb{R} \times_{f(r)} S^1$
where
\be
f(r)=2Arctan(-r)+\pi.
\ee
Since $\lim_{r\to\infty}f(r)=0$, $\pi_{slip}(M)=\pi_1(M)$
and the covering spectrum is empty.

Now let $(X_i,x_i)=(M,p_i)$ where $r(p_i)=r_i \to -\infty$.
Note that $\bar{B}_{x_i}(R)$ is then equipped with a warped product metric
and
\be
f_i(r)=f(r-r_i)=2Arctan(-r+r_i) +\pi
\ee
which converges uniformly on $[-R,R]$ to
\be
f_{\infty}(r)=2\pi.
\ee
Thus the pointed Gromov Hausdorff limit is the standard cylinder
whose covering spectrum is $\{\pi\}$.  As above the
cut-off covering spectra of these examples is empty
both for the $X_i$ and the limit space.
\end{example}

Next we construct an example where an element of the covering spectrum
disappears in the limit without decreasing to 0.
This issue is not immediately solved
by using the cut-off covering spectrum.

\begin{example}  \label{disappear}
Let $M^2$ be a
cylinder with a small handle near a point $p$. Let
$(X_i,x_i)=(M,p_i)$ where $d(p_i,p)\to \infty$. Then
$CovSpec(X_i)=CovSpec(M)$ since the covering spectrum does not
depend on the base point and the spectrum has its first element,
$\lambda_1<\pi$ corresponding to the small handle.
 Yet for all $R>0$
just take $N_R$ large enough that
\be
d(p_i,p) \ge 2R \ \ \ \  \forall i\ge N_R.
\ee
Then $\bar{B}_{p_i}(R)$ are all isometric to balls of radius
$R$ in a cylinder.  So $(X_i,x_i)$ converge to a cylinder
with a point.  So the covering spectrum of the limit space
does not include $\lambda_1$ and its only element is $\pi$.
So we have locally compact $X_i$ converging to locally
compact $X$ with \be \delta_i \in CovSpec(X_i) \textrm{ such that }
\delta_i=\lambda_1 \rightarrow
\delta \notin CovSpec(X). \ee
In fact we have
\be \delta_i = \lambda_1 \in CovSpec_{cut}(X_i) \textrm{ such that }
\delta_i
\rightarrow \delta \notin CovSpec_{cut}(X). \ee

On the other hand any $R>0$,
there exists $N_R>0$ such that
$CovSpec_{cut}^R(X_i,x_i)=\{\pi\}$ for all $i\ge N_R$
because the handle is
located outside $\bar{B}_{x_i}(R)$.
\end{example}

Finally we have the possibility that elements of the covering
spectrum can grow to infinity.  In this example we see that in
essence a hole could expand until it snaps and is no longer a hole
in the limit space:

\begin{example}\label{expandtoinfinity} \label{ex5.5}
Let $(X_r,x_r)$ be formed where $X_r$ is a unit interval
with $x_r$ on one end and a circle of circumference $2\pi r$
attached to the other end with a half line attached on the
opposite side of the circle.   Then $CovSpec(X_r)=\{\pi r \}$.

Note that if one takes a sequence of $r_i$ diverging to
infinity, $(X_{r_i},x_{r_i})$, converges in the pointed
Gromov-Hausdorff sense to $(X_\infty,x_\infty)$ where
$X_\infty$ is a unit interval attached to $x_\infty$ at
one end and two half lines at the other end.  So $X_\infty$
is simply connected and has an empty covering spectra.

This example is not
simplified by using the cut-off covering spectra.
In fact for any $R\ge 1$
$CovSpec_{cut}^R(X_r)=\{\pi r\}$ and so $CovSpec_{cut}(X_r)=\{\pi r\}$.
\end{example}

\subsection{Convergence of the $R$ cut-off covering spectrum}\label{sect5.2}

In light of the above examples, it is natural to try to to prove
continuity of the $R$ cut-off covering spectra and then perhaps to
apply this continuity to prove some form of continuity for the
cut-off covering spectrum.  Surprisingly the statement of the
continuity theorem for the cut-off spectrum is somewhat tricky:

\begin{theo} \label{covcutconv} \label{thm5.2}
Let $(X_i,x_i)$ be complete locally compact length spaces 
converging in the pointed
Gromov-Hausdorff sense to a locally compact space $(X,x)$.
Bounded elements do not disappear:
if we have a converging sequence,
\be\label{ccca}
\delta_i \in CovSpec^{R_1}_{cut}(X_i,x_i)
\textrm{ and } \delta_i \rightarrow \delta >0, \textrm{ then }
\delta \in CovSpec^{R_1}_{cut}(X,x).
\ee
Nor do elements suddenly appear: for any $R_2>R_1$ and if we have an
element,
\be\label{cccb}
\delta \in CovSpec^{R_1}_{cut}(X,x)
\textrm{ there are }
\delta_i \in CovSpec^{R_2}_{cut}(X_i,x_i)
\textrm{ such that } \delta_i \rightarrow \delta.
\ee
\end{theo}

Examples~\ref{shrinkto0} and~\ref{expandtoinfinity}
demonstrate why one must assume $\delta_i$ converge in $(0,\infty)$.
We now present examples demonstrating why we cannot take
$R_1=R_2$ in (\ref{cccb}).

\begin{example} \label{10.3}
Let $(X_r,x_r)$ be formed by attaching a line segment
of length $r$ to a circle of circumference $2\pi$, and
then continuing with a half line on the opposite side of the circle.
The point $x_r$ will be the endpoint of the line segment
not attached to the circle.  If $r_i \to r_\infty$ it is
easy to see that $(X_{r_i},x_{r_i})$ converges
to $(X_{r_\infty},x_{r_\infty})$.

Note that $CovSpec(X_{r})=\{\pi\}$ and so does $CovSpec_{cut}(X_r)$.
On the other hand, $CovSpec_{cut}^R(X_r)=\emptyset$ when $r>R$ because
then the circle is contained in $X_r \setminus \bar{B}_{x_r}(R)$.
Otherwise $CovSpec_{cut}^R(X_r)=\{\pi\}$.  Thus
the sequence $X_{r_j}$ with $r_j$ decreasing to $R_1$
has
\be \label{10.1aequ}
\delta=\pi\in CovSpec_{cut}^{R_1}(X_{r_\infty})
\ee
but $CovSpec_{cut}^{R_1}(X_{r_i})=\emptyset$.  However, taking
$R_2>R_1$, eventually we have $r_i<R_2$ so we have
\be \label{10.1aequ2}
\delta_i=\pi\in CovSpec_{cut}^{R_2}(X_{r_i}).
\ee
\end{example}

The next example also illustrates the same phenomenon
with a distinct cause:

\begin{example} \label{10.3a}
Let $M$ be a warped product manifold of the form
$\Bbb{R} \times _f S^1$ where $f(t)=e^{-t^2}$.  Fix
$p$ in the level $t=0$.

Let $X_r = \bar{B}_p(r)$.
So it is a closed ball and if we wish to make it noncompact, we just
attach a half line to it.  We give
it the induced length metric from $M$.

Let $r_i$ decrease to some $r_\infty > \pi$.
Then $X_i=X_{r_i}$ converges
to $X_\infty=X_{r_\infty}$.

Let $R_1 =r_\infty$.  The $R_1$-cutoff covering spectrum of $X_\infty$
includes $\delta$ equal to half the length of one of the
components of the boundary of $B(p,R_1)$, because this curve
is not homotopic to anything outside $\bar{B}_{p}(R_1)$.  However the
$R_1$-cutoff covering spectra of the $X_i$ are all empty because
the loop is homotopic to a loop in $\partial B_p(r_i)$ which
is outside $B(p,R_1)$.  So once again we need $R_2>R_1$ and need
to wait for $r_i<R_2$ to get the cut-off covering spectra to
converge.

One might also construct manifolds $M_i$ converging to $X_\infty$
by taking smoothed tubular neighborhoods of the $X_i$ in
five dimensional Euclidean space.
\end{example}

Note that in the above examples, if one were to take $r_j$ increasing
to $R_0$ the covering spectrum are all $\{\pi\}$.  The difficulty
arises because $r_j$ decreasing to $R_\infty$ are leaving the open
set $(R_\infty,\infty)$ in the limit.

At first we thought we needed to take $R_2>R_1$ in (\ref{ccca}) as
well as (\ref{cccb}) but due to the lack of examples proving this
was necessary, we investigated further and discovered we could boost
our proof of (\ref{ccca}) using the local compactness of the limit space.

In order to prove this theorem we need extend several results for covering
spaces of compact spaces to $R$-cutoff spaces.  The first
is an adaption of Theorem 3.4 in \cite{SoWei1}.

\begin{prop}  \label{theodelsurj}
Let $B(p_i,s_i) \subset B(p_i,S_i)\subset Y_i, i=1,2$ be
 balls each with induced length metrics.  Let
$G(p_1,s_1,S_1,\delta_1)$ be the group of
deck transformations of
$\tilde{B}(p_1, S_1)_{cut}^{\delta_1, s_1}$.

If there is a
pointed $\epsilon$-Hausdorff approximation
$f:B(p_1,S_1) \ra B(p_2,S_2)$
 then
for any $\delta_1> 10 \epsilon$ and $\delta_2>\delta_1+10\epsilon$
and $s_2 < s_1-5\epsilon$,
there is a surjective homomorphism,
\be
\Phi: G(p_1,s_1,S_1,\delta_1) \to G(p_2,s_2,S_2,\delta_2).
\ee
\end{prop}

\noindent {\bf Proof of Proposition~\ref{theodelsurj}}:
We begin by describing a map for closed curves.
For a closed curve $\gamma: [0,1] \rightarrow B(p_1,S_1)$ with $\gamma(0)=
\gamma(1) =p_1$, construct a $5\epsilon$-partition of $\gamma$ as follows.
On
$\Gamma := \gamma ([0,1])$
choose a partition $0 =t_0 \leq t_1 \leq \cdots
\leq t_m = 1$ such that for $x_i = \gamma(t_i)$, one has
$d(x_i, x_{i+1}) < 5\epsilon$ for $i= 0, \cdots, m-1$.
$\{x_0,
\cdots, x_m\}$ is called a $5\epsilon$-partition of $\gamma$.

Let  $y_m = y_0=p_2$ and for each $x_i$, we set
$y_i =f(x_i), i=1,\cdots,m-1$. Connect $y_i$ and $y_{i+1}$ by minimal
geodesics in $B_{p_2}(S_2)$.  This yields a closed curve
$\bar{\gamma}$ in $B(p_2,S_2)$ based at $p_2$ consisting of $m$ minimizing
segments each having length $\le 6\epsilon$.

Any $\alpha \in G(p_1,s_1,S_1,\delta_1)$
can be represented by some rectifiable closed curve
$\gamma$ in $B(p_1,S_1)$, so we can hope to define
\[
\Phi (\alpha) = \Phi ([\gamma]) := [\bar{\gamma}] \in
G(p_2,s_2,S_2,\delta_2).
\]
First we need to verify that $\Phi$ doesn't depend on
the choice of $\gamma$ such that $[\gamma]=\alpha$.

Using the facts that $18\epsilon<\delta_2$ and loops which fit
in balls of radius $\delta_2$ do not effect the representative
of a class in $G(p_2,s_2, S_2,\delta_2)$,
one easily see that $[\bar{\gamma}]$ doesn't depend on the choice of
minimizing curves $\bar{\gamma}_i$,
nor on the special partition $\{x_1, \cdots, x_m \}$ of $\gamma ([0,1])$.

Moreover using additionally the uniform continuity of a homotopy one
can similarly check that if $\gamma$ and $\gamma'$ are homotopic in
$B(p_1,S_1)$, then $[\bar{\gamma}]=[\bar{\gamma}']$ in
$G(p_2,s_1-5\epsilon,S_2,\delta_2)$. That is, we can take a homotopy
$h:[0,1]\times[0,1]\to B(p_1,S_1)$, we can take a grid on
$[0,1]\times[0,1]$ small enough that  homotopy maps the grid points
to points $x_{i,j}$ that are less than $5\epsilon$ apart from the
images of their grid neighbors.  Then we take $y_{i,j}=f(x_{i,j})$
and connect neighbors according to the rules in the first paragraph.
Finally we use the argument in the paragraph above this to see that
the net created using the $y_{i,j}$ is a $\delta_2$ homotopy so
$[\bar{\gamma}]=[\bar{\gamma}']$ in
$G(p_2,s_2-5\epsilon,S_2,\delta_2)$. Thus we see that $\Phi$ is a
homomorphism from $\pi_1(B(p_1,S_1), p_1)$ to $G(p_2,s_2,S_2,\delta_2)$.
However $\alpha \in
G(p_1,s_1,S_1,\delta_1)$ not $\pi_1(B(p_1,s_1), p_1)$.

Suppose $\gamma_1$ and $\gamma_2$ are both representatives of $\alpha
\in G(p_1,s_1,S_1,\delta_1)$.  Then $\gamma_1 *\gamma_2^{-1}$ is, in
$B(p_1,S_1)$,
homotopic to a loop $\gamma_3$ generated by loops of the
form  $\alpha * \beta * \alpha^{-1}$, where $\beta$ is a closed path lying
in a ball of radius $\delta_1$ or
in $B(p_1,S_1)\setminus \bar{B}(p_1, s_1)$.
So $[\bar{\gamma_1}]=[\bar{\gamma_3}] * [\bar{\gamma_2}]$
in $\pi_1(B(p_1,s_1), p_1)$.  So
we need only show that $[\bar{\gamma_3}]$ is trivial in
$G(p_2,s_2,S_2,\delta_2)$.

In fact $\bar{\gamma_3}$ can be chosen  as follows.
The $y_i$'s corresponding to the $x_i$'s from the $\beta$ segments of
$\gamma_3$  are all within $\delta_1+\epsilon$ of a common point
and the minimal geodesics between them are within
$\delta_1+(1 +6/2) \epsilon < \delta_2$.  Furthermore, the  $y_i$'s
corresponding to the $x_i$'s from the $\alpha$ and $\alpha^{-1}$ segments
of the curve can be chosen to correspond.  Thus  $\bar{\gamma_3}$
is generated by loops of the
form  $\alpha * \beta * \alpha^{-1}$ lying in $B(p_2,S_2)$, where $\beta$ is
a
closed path lying
in a ball of radius $\delta_2$ or $B(p_1,S_2)\setminus B(p_1,
s_1-5\epsilon)$
and $\alpha$ is a path from $p_2$ to $\beta(0)$.
So it is trivial.

Last, we need to show that $\Phi$ is onto. If
$\bar{\alpha} \in G(p_2,s_2,S_2,\delta_2)$,
it can be represented by some
rectifiable closed curve $\sigma$ in $B(p_2,S_2)$ based at $p_2$.  Choose an
$\epsilon$-partition
$\{y_0, \cdots, y_m\}$ of $\sigma$. Since
$f: B(p_1,S_1) \ra B(p_2,S_2)$ is an $\epsilon$-Hausdorff
approximation, there are $x_i \in B(p_1,s_1)$, $y_i'=f(x_i) \in B(p_2,S_2)$
where $y_0' = y_m' =p_2$, $x_0 = x_m =p_1$ and
$d_{B(p_2,S_2)}(y_i,y_i') \le \epsilon$.
Connect $y_i', y_{i+1}'$ with a length minimizing curve in $B(p_2,S_2)$;
this
yields a piecewise length minimizing closed curve $\sigma'$ in $B(p_2,S_2)$
based at $p_2$, each segment has length $\le 3\epsilon$. So $[\sigma'] =
[\sigma]$ in $G(p_2,s_1-5\epsilon,S_2,\delta_2)$. Now connect
$x_i, x_{i+1}$ by length minimizing
curves in $B(p_1,S_1)$ this yields a piecewise length minimizing
$\gamma: [0,1] \rightarrow B(p_1, S_1)$ with base point $p_1$, each segment
has
length $\le 4\epsilon$. So the curve $\gamma$ allows a
$5\epsilon$-partition and $[{\gamma}]\in G(p_1, s_1,S_1,\delta_1)$.
By the construction, $\Phi([\gamma])=\bar{\alpha}$.

Therefore $\Phi$ is surjective.
\qed

\begin{prop}  \label{delta-precom}
 If a sequence of complete locally compact length spaces
$X_i$ converges to a  length
 space $X$ in the Gromov-Hausdorff topology, then for any
$\delta >0, R>0, r>3R$ there is a
 subsequence of $X_i$ and a sequence $r_i \to r$ such that
$\tilde{B}(x_i, r_i)_{cut}^{\delta, R}$ also converges in the
 pointed Gromov-Hausdorff topology.
Moreover,
the limit space $B(x,r)_{cut}^{\delta, R}$ is a
covering space of $B(x,r)$ satisfying
 \be
 \tilde{B}(x,r)_{cut}^{\delta, R}
\rightarrow B(x,r)_{cut}^{\delta, R} \rightarrow
\tilde{B}(x,r)_{cut}^{\delta', R'}
 \ee
 for all $0< R' < R$ and $\delta'> \delta$.
 \end{prop}

\Pf
By the Appendix of \cite{SoWei2} we know that for a sequence
$r_i$ converging to $r$, $B(x_i,r_i)$ converge with the induced length
metric to $B(x,r)$.

By \cite{SoWei3}[Proposition 7.3] and the fact that the closed balls
$B(x_i,r_i)$ are compact sets, we know that $\tilde{B}(x_i,r_i)^\delta$
have a converging subsequence.  So by Gromov's compactness theorem,
they have a uniform bound $N(a,b)$, the number of disjoint
balls of radius $a$ in a ball of radius $b$.   By
Proposition~\ref{cutcoverexists},
$\tilde{B}(x_i, r_i)^{\delta}_{cut}$ covers
$\tilde{B}(x_i, r_i)_{cut}^{\delta, R}$, so
$N(a,b)$ can be used to count balls in
$\tilde{B}(x_i, r_i)_{cut}^{\delta, R}$ as well.  So by Gromov's
compactness theorem, a subsequence of these spaces converges and
we will denote the limit space:
${B}(x, r)_{cut}^{\delta, R}$.

To complete the proof we adapt Theorem 3.6 of \cite{SoWei1}.
The fact $\tilde{B}(x_i, r_i)_{cut}^{\delta, R}$ were isometries
on balls of radius $\delta$ and outside $\bar{B}_p(R)$
guarantees that the limit is as well,
so ${B}(x, r)_{cut}^{\delta, R}$ is a covering space for
$B(x,r)$.

The isometries also guarantee it is covered by
$\tilde{B}(x,r)_{cut}^{\delta,R}$.  This can be seen using the
Unique Lifting Theorem (c.f. \cite{Massey} Lemma 3.1, p123)
and noting that if $C$ is a closed curve in $B(x,r)$ whose lift to
$\tilde{B}(x,r)_{cut}^{\delta,R}$
then it is homotopic to a curve which is created from curves
of the form $\alpha\cdot\beta\cdot \alpha^{-1}$ where the
$\beta$ are either in a ball of radius $\delta$ or outside $\bar{B}_x(R)$.
So its lift to $B(x,r)_{cut}^{\delta, R}$ is also closed since $\pi^\delta$
is an
isometry on $\delta$-balls and an isometry outside $\bar{B}_x(R)$.
Therefore $\tilde{B}(x,r)_{cut}^{\delta,R}$ covers ${B(x,r)}^{\delta,R}$

To complete the proof we apply the Unique Lifting Theorem
by contradiction.
We assume there is $\delta'>\delta$ and $R'<R$ and $C$ is a curve
which lifts closed to $B(x,r)_{cut}^{\delta, R}$ but lifts open to
$\tilde{B}(x,r)_{cut}^{\delta', R'}$.
Since this lift of $C$ is not closed,
$[C] \in G(x,r,R',\delta')$ is nontrivial.

Let $\epsilon>0$ be chosen sufficiently small that \be
\epsilon<\min\{\delta/10,(\delta-\delta')/10,(R-R')/5\}. \ee Take
$i$ sufficiently large that we have an $\epsilon$-Hausdorff
approximation $f_i:B(x_i,r_i) \ra B(x,r)$.  Applying
Proposition~\ref{theodelsurj}, we know there are a surjective
homomorphisms, $\Phi: G(x_i,r_i,R,\delta) \to G(x,r,R',\delta')$, so
there are closed loops $C_i\in B(x_i,r_i)$ such that
$\Phi([C_i])=[C]$.

By the construction of $\Phi$, $C_i$ can be chosen so
these lifted curves $\tilde{C_i}$ converge to the
lift of the limit of the curves, $\tilde{C}$ in $B(x,r)_{cut}^{\delta,R}$
and
\be \label{0to1}
d_{B(x,r)_{cut}^{\delta,R}}(\tilde{C}(0), \tilde{C}(1))=
\lim_{i\to\infty}d(\tilde{C}_i(0), \tilde{C}_i(1)).
\ee
However the $[C_i]$ are nontrivial, so their lifts to
$\tilde{B}(x_i,r_i)^{\delta,R}$
run between points $\tilde{C}_i(0)\neq\tilde{C}_i(1)$ satisfying
\be
d(\tilde{C}_i(0), \tilde{C}_i(1))\ge \delta.
\ee
Combining this with (\ref{0to1}), we see that
$\tilde{C}$ is not closed and we have a contradiction.
\qed

At this point we could imitate the proof of Theorem 8.4 in \cite{SoWei3}
to prove Theorem~\ref{covcutconv} for $X_i$ which are compact balls.
However, this would not help us prove Theorem~\ref{covcutconv} for
noncompact spaces as the cut-off covering spectrum of a ball does
not match the cut-off covering spectrum of the space. Recall
Examples~\ref{expandtoinfinity} and~\ref{cylinder} demonstrate that
not only can holes become increasingly large, but homotopies may as
well. One needs to control such phenomenon to complete the proof.

\vspace{.2cm} 

\noindent{\bf Proof of Theorem~\ref{covcutconv}}:
In order to prove the first statement (\ref{ccca})
we first prove that given any $R_2>R_1$
if
\be\label{cccc}
\delta_i \in CovSpec^{R_1}_{cut}(X_i)
\textrm{ and } \delta_i \rightarrow \delta >0, \textrm{ then }
\delta \in CovSpec^{R_2}_{cut}(X).
\ee
Later we will boost this result to (\ref{ccca}).

Assume
\be
\delta_i \in CovSpec^{R_1}_{cut}(X_i)
\ee
and $\delta_i \rightarrow \delta >0$. By Lemma~\ref{globaltolocal},
$\delta_i \in CovSpec_{cut}^{R_1}(B(x_i,r))$ for $r \ge 3(R_1+2\delta_i)$.
So $\tilde{B}(x_i,r)_{cut}^{\delta_i,R_1} \to
 \tilde{B}(x_i,r)_{cut}^{\delta',R_1}$ is nontrivial for all
$\delta'>\delta_i$.  So for all $\delta'>\delta>0$ and
 $\epsilon\in (0, \delta)$
 we have $\delta-\epsilon<\delta_i<\delta'$  for $i$ sufficiently large and
 \be
\tilde{B}(x_i,r)_{cut}^{\delta-\epsilon,R_1} \to
\tilde{B}(x_i,r)_{cut}^{\delta',R_1}
\ee
is nontrivial.  Now take the limit as $i \to \infty$ and we get
\be
B(x,r)_{cut}^{\delta-\epsilon,R_1} \to
B(x,r)_{cut}^{\delta',R_1}
\ee
is nontrivial.

 This is true for all $\epsilon\in (0,\delta)$ and $\delta'>\delta$.
 By the properties of limit covers in Proposition~\ref{delta-precom}
we have for all $\epsilon\in (0,\delta)$, $\delta''>\delta'$, 
and $R'\in(R_1,R_2)$,
\be
\tilde{B}(x,r)_{cut}^{\delta-\epsilon,R_1} \to
B(x,r)_{cut}^{\delta-\epsilon,R_1} \textrm{ and }
B(x,r)_{cut}^{\delta',R_1}
 \to
\tilde{B}(x,r)_{cut}^{\delta'',R'}.
\ee
Therefore $\tilde{B}(x,r)_{cut}^{\delta-\epsilon,R_1} \to
\tilde{B}(x,r)_{cut}^{\delta'',R'}$
is nontrivial.

By Proposition~\ref{iflex}, we then know
that since $B(x,r)$ is compact and $R_2>R'$ we have
\be
CovSpec^{R_2}_{cut}(B(x,r))\cap [\delta-\epsilon,\delta'')\neq \emptyset.
\ee
Taking $\epsilon$ to $0$ and $\delta''$ to $\delta$, we get
\be
\delta \in
CovSpec^{R_2}_{cut}(B(x,r)
\ee
This is true for all sufficiently large $r$, so
by Lemma~\ref{uniflocaltoglobal},
\be
\delta \in CovSpec^{R_2}_{cut}(X)
\ee
which completes proof of (\ref{cccc}).

We now boost the statement (\ref{cccc}) to prove (\ref{ccca}).
Again fix $R_1 >0$. Suppose
\be
\delta_i \in CovSpec^{R_1}_{cut}(X_i)
 \ee
 and $\delta_i \rightarrow \delta >0$,
 Let $X$ be the Gromov-Hausdorff limit of the $X_i$.

Proposition~\ref{iflex} says that that for $R_2=R_i$ sufficiently
 close to $R_1$
 \be
 \tilde{X}^{\delta,R_2}_{cut}= \tilde{X}^{\delta,R_1}_{cut}
 \ee
 putting this together with Proposition~\ref{Rflex} says
 \be
 \delta \notin CovSpec^{R_2}_{cut}(X) \setminus CovSpec^{R_1}_{cut}(X)
 \ee
 We apply (\ref{cccc}) to say
 \be
 \delta \in CovSpec^{R_2}_{cut}(X).
 \ee
 But then
 \be
 \delta \in CovSpec^{R_1}_{cut}(X),
 \ee
 which gives us (\ref{ccca}).

Now we prove the second statement (\ref{cccb}):
given
\be \label{setdelta}
\delta \in CovSpec^{R_1}_{cut}(X)
\ee
and any $R_2>R_1$, show there exists
\be
\delta_i \in CovSpec^{R_2}_{cut}(X_i)
\ee
such that $\delta_i \rightarrow \delta$.

We assume on the contrary that there is a gap:
\be\label{gap}
\exists \epsilon>0
\textrm{ such that }
CovSpec^{R_2}_{cut}(X_i) \cap (\delta-2\epsilon,\delta+2\epsilon)=\emptyset.
\ee
By Lemma~\ref{uniflocaltoglobal}, for $r\ge 3(R_2 +2 \delta +4\epsilon)$,
\be
CovSpec^{R_2}_{cut}(B(x_i,r)) \cap
(\delta-2\epsilon,\delta+2\epsilon)=\emptyset.
\ee

By Lemma~\ref{subsets} we then have for any $R_2' \le R_2$:
\be
CovSpec^{R'_2}_{cut}(B(x_i,r)) \cap
(\delta-2\epsilon,\delta+2\epsilon)=\emptyset.
\ee

So the covering
\be
\tilde{B}(x_i,r)_{cut}^{\delta-\epsilon, R_2'}
\to \tilde{B}(x_i,r)_{cut}^{\delta+\epsilon, R_2'}
\ee
 is trivial. By Proposition~\ref{delta-precom}
we have a subsequence of the $i$ such that:
\be
\tilde{B}(x_i,r)_{cut}^{\delta-\epsilon, R_2'} \to
B(x,r)_{cut}^{\delta-\epsilon, R_2'}
\ee
and
\be
\tilde{B}(x_i,r)_{cut}^{\delta+\epsilon, R_2'} \to
B(x,r)_{cut}^{\delta+\epsilon, R_2'}.
\ee
since the sequence of the covering map is trivial,  the covering limit map
\be  \label{applythis}
B(x,r)_{cut}^{\delta-\epsilon, R_2'}  \to B(x,r)_{cut}^{\delta+\epsilon,
R_2'}
\ee
 is also trivial.

By Proposition~\ref{delta-precom} for any $R_2 \ge R_2' >R_1$
 \be \label{applythis2}
B(x,r)_{cut}^{\delta-\epsilon, R_2'}
 \rightarrow  \tilde{B}(x,r)_{cut}^{\delta, R_1}
\rightarrow  \tilde{B}(x,r)_{cut}^{\delta+\epsilon, R_1}
\rightarrow B(x,r)_{cut}^{\delta+\epsilon, R_1}.
 \ee

By Proposition~\ref{iflex},  for any $R_2'>R_1$ sufficiently close to $R_1$,
the covering
\be
\tilde{B}(x,r)_{cut}^{\delta- \epsilon, R_2'} \to
\tilde{B}(x,r)_{cut}^{\delta+2\epsilon, R_1}
\ee
is trivial.
Using Proposition~\ref{delta-precom}, we have for $R_2'' <R_2'$, the
covering
\be
B(x,r)_{cut}^{\delta-\epsilon, R_2''} \ra B(x,r)_{cut}^{\delta+\epsilon,
R_1}
\ee
is trivial.

Apply this $R_2''$ to (\ref{applythis2}) we get
 trivial covers in {\ref{applythis2}).
So
$\delta \notin CovSpec^{R_1}_{cut}(B(x,r))$.

By Lemma~\ref{globaltolocal},
$\delta \notin CovSpec^{R_1}_{cut}(X)$. That is a contradiction.
\qed

\subsection{Convergence of the cut-off covering spectrum} \label{sect5.3}

Theorem~\ref{covcutconv} combined with Proposition~\ref{bigcupconverse}
gives
the following result that elements in the cut-off covering spectrum
do not suddenly appear in limits.  Example~\ref{disappear}
demonstrates that elements
of the cut-off covering spectrum can disappear in the limit by sliding
out to infinity.  Unlike the $R$ cut-off covering spectrum, all handles are
now visible.

\begin{theo} \label{noappear}\label{thm-no-appear}
Let $(X_i,x_i)$ be complete locally compact length spaces
converging in the pointed
Gromov-Hausdorff sense to a locally compact space $(X,x)$, then
\be
\textrm{ for any }
\delta \in CovSpec_{cut}(X), \textrm{ there is }
\delta_i \in CovSpec_{cut}(X_i)
\ee
 such that $\delta_i \rightarrow \delta$.
\end{theo}

This provides an immediate application:

\begin{coro}\label{limsimp}
If $X_i$ are simply connected locally compact length spaces
converging in the pointed
Gromov-Hausdorff sense to a locally compact space $(X,x)$
then $CovSpec_{cut}(X)=\emptyset$.
\end{coro}

\noindent{\bf Proof of Theorem~\ref{thm-no-appear}:}
 If  $\delta \in CovSpec_{cut}(X)$, by Proposition~\ref{bigcupconverse},
\be
\delta \in Cl_{lower} \cup_{R>0} CovSpec^R_{cut}(X).
\ee
So there are $R_k$ increasing to infinity and
\be
\delta_k \in CovSpec^{R_k}_{cut}(X)
\ee
 such that $\delta_k \rightarrow \delta$ .
By Proposition~\ref{covcutconv}
and $R_{k+1}>R_k$, for each $\delta_k$, we have
\be
\delta_k^i \in CovSpec^{R_{k+1}}_{cut}(X_i) \subset CovSpec_{cut}(X_i)
\ee
such that $\delta_k^i  \rightarrow \delta_k$.

By a diagonal process, we have
\be
\delta_i=\delta_{k_i}^i \in CovSpec_{cut}(X_i)
\ee
such that $\delta_i \rightarrow \delta$.
\qed

\begin{ques}\label{questproperconv}
Is local compactness a necessary condition
in our convergence theorems [Theorem~\ref{thm5.2} and
Theorem~\ref{thm-no-appear}]?
This condition is used in a few crucial
steps of the proof.  It is used in Lemma~\ref{lem3.7} to apply the
pigeonhole
principle to control the lengths of shortest representative curves.
It is also used in Proposition~\ref{delta-precom} to prove that the delta
covers of balls converge in the Gromov-Hausdorff sense when these
balls converge in the Gromov-Hausdorff sense.    This proposition
is based on a result in \cite{SoWei3} which requires compactness.
Proposition~\ref{iflex} which requires compactness is applied to balls in
the proof.  Finding examples demonstrating the necessity of local
compactness or compactness in any of these results would be
of interest.
\end{ques}

Recall the pulled ribbon construction introduced in Section~\ref{sectribbon}
was used to obtain important examples which are not locally compact.
Here however the ribbon construction does not immediately help:

\begin{rmk}\label{rmkproperconv}
Suppose one were to attach pulled ribbons to the sequence
of manifolds in Example~\ref{ex5.3} in an attempt to prove that
local compactness is necessary in Theorem~\ref{noappear}.
The difficulty is that as soon as
the spaces are bounded, there is no way to effectively use a base point
to differentiate the spaces from one another.
Thus the sequence is just a repeating space and converges to itself.  It
does not produce a counter example.  The same effect happens if we
try to attach pulled ribbons to capped cylinders which converge
to cylinders in the pointed Gromov Hausdorff sense (see Example~\ref{ex5.2}).

If we choose $M_j$ to be warped products with two cusps that are
isometric to cylinders on $[-j,j]$ such spaces would converge in the
pointed Gromov-Hausdorff sense to a cylinder, but when we attach pulled
ribbons to them, the sequence does not even converge: each space is
a definite Gromov-Hausdorff distance apart from each other.  The
only reason the $M_j$ converged was because the pointed GH convergence
only saw the center cylindrical region, but when the whole space is
bounded, the whole space needs to behave in a uniform way.
\end{rmk}

\subsection{Applications of Convergence} \label{sect5.4}

In this section we observe the following topological
consequence of our convergence results:

\begin{theo}\label{thm5.7}
If $X_i$ are complete locally compact  length spaces that satisfy the loops to
infinity property and converge in the pointed Gromov-Hausdorff
sense to a locally compact  semi-locally simply connected limit space $X$
then either $X$ has at least two ends
or $X$ has the loops to infinity property.
\end{theo}

\Pf
By Theorem~\ref{thmloopstoinfty}, $CovSpec_{cut}(X_i)$ are trivial. 
So, applying
Theorem~\ref{thm-no-appear}, we see $CovSpec_{cut}(X)$ must be trivial.
To complete the proof we just apply Theorem~\ref{thm-one-end}.
\qed

One can think of this theorem as the complete version of
the theorem in \cite{SoWei1} which says that compact Gromov-Hausdorff
limits of simply connected compact manifolds are simply connected.

\begin{example} \label{simplyconnected}
The \cite{SoWei1} theorem is not true for noncompact limits with pointed
Gromov-Hausdorff convergence as can be seen by taking sequences
of ellipsoids $M_j^2$ which stretch out to a cylinder $S^1\times\Bbb{R}$
or $M_j^2\times\Bbb{R}$ converging to $S^1\times\Bbb{R}^2$.
Thanks to our new theorem we see that while holes may form in a limit
they cannot be handles.
\end{example}

\begin{example}\label{GHtwoends}
Notice if one takes a disk and stretches two points out to infinity
then the limit is a disk with two cusps, which is no longer simply
connected.  Nor does it has the loops to infinity property.  This is
because a loop wrapping once around each cusp is not homotopic to
loops approaching infinity.  However the fundamental group of the
space is generated by elements with the loops to infinity property.
\end{example}

\begin{example}\label{pi2}
Note that one can have compact $M_j$ with $\pi_2(M_j)$ nontrivial converging to
a space with nontrivial $\pi_2$.  This can be seen by taking $M_j$
diffeomorphic to the plane with warped product metrics \be dr^2 +
f_j^2(r) d\theta^2 \ee where $f(r)=r ((1-r)^2 +(1/k) )$, so that the
Gromov-Hausdorff limit as $k\to\infty$ is homoemorphic to a sphere
attached to a plane. So we cannot hope to control higher homotopy,
although an investigation of \cite{ShSo1} reveals a close
relationship between the loops to infinity property and the
codimension one integer homology of the space.
\end{example}

\subsection{Tangent Cones at Infinity}\label{sect5.5}

A complete noncompact space, $X$, is said to have a tangent cone at
infinity if the Gromov-Hausdorff limit of a sequence of inward rescalings
$(X/{r_j},x)$ with $r_j \to \infty$, has a limit in the
pointed Gromov-Hausdorff sense.  While this limit space is called a cone,
it is not a metric cone except in very special situations, like when
$X$ has nonnegative sectional curvature \cite{BBI}.  In fact the tangent
cone at infinity of a manifold need not even be simply connected
as can be seen in this well-known example:

\begin{example}
Let $M^2$ be created by taking a cone, smoothing off the tip
and adding handles, $r_iH$, at a distance $r_i$ from the old tip.  We
write $r_iH$ because we are rescaling the handle $H$ by $r_i$, so
that the handles are growing.  Then $M^2/r_i$ converges to a cone
with a handle attached at a distance $1$ from the tip.  If
$\lim r_{i+1}/r_i =\infty$, then the tangent cone has only one
handle,  but if $\lim r_{i+1}/r_i= d$, then the tangent cone
has infinitely many handles located at $\{d^j: j\in \Bbb{Z}\}$, so
the tangent cone at infinity, $Y$, has
locally infinite topological type at its tip.
Furthemore $Y$ has no universal cover
and $CovSpec(Y)=CovSpec_{cut}(Y)$ have infinitely many elements.
\end{example}

\begin{rmk}
Menguy has created similar examples demonstrating that the
tangent cone at infinity of a manifold with nonnegative
Ricci curvature can have locally
infinite topological type, although his
examples are simply connected because his handles are
higher dimensional (c.f. \cite{ShSo2}).  In \cite{SoWei3}
we proved the tangent cones at infinity of manifolds with
$Ricci \ge 0$ have universal covers.
\end{rmk}

Using our results we can prove

\begin{theo} \label{tannoappear}
If $X$ is a complete locally compact length space and
$CovSpec_{cut}(X)$ is bounded
then any tangent cone at infinity for
$X$ has a trivial cut-off covering spectrum.
\end{theo}

First note that the following lemma holds

\begin{lemma}\label{rescale}
If we rescale a metric space $X$ to get a new
metric space $X/r$
then the elements of the covering spectrum
and cut-off covering spectrum scale
proportional to the distance:
\be
CovSpec(X/r) = CovSpec(X)/r
\ee
and
\be
CovSpec_{cut}(X/r) = CovSpec_{cut}(X)/r.
\ee
Furthermore, we have
\be
CovSpec_{cut}^{R/r}(X/r) = CovSpec_{cut}^{R}(X)/r.
\ee
\end{lemma}

This lemma follows immediately from the definitions.

\noindent{\bf Proof of Theorem~\ref{tannoappear}:}
If we rescale a space dividing the metric
by $r$ then by Lemma~\ref{rescale} we have
\be
CovSpec_{cut}^{R/r}(X/r) = CovSpec_{cut}^{R}(X)/r
\ee
So applying Proposition~\ref{bigcupconverse} we have
\be
CovSpec_{cut}(X/r) = CovSpec_{cut}(X/r) \subset [0, Max(CovSpec(X))/r]
\ee
Any tangent cone at infinity, $Y$, is the Gromov-Hausdorff limit of $X/r_i$
with $r_i\to\infty$, so by Theorem~\ref{thm-no-appear},
$
CovSpec_{cut}(Y) \subset \{0\}
 $
and is, thus, trivial.
\qed

This proposition implies that the tangent cones at infinity
of manifolds with bounded covering spectra have trivial
cutoff covering spectra.  However, they need not have trivial
covering spectra even when the manifold has a trivial covering
spectrum:

\begin{example}
Let $M$ be the length space constructed by attaching a sequence
of widening cylinders to a plane as follows.  Take a flat Euclidean
plane and remove disks of radius $2^j/4$ about the points
$(2^j,0)$ where $j=1,2,3...$.  Now attach standard cylinders of radius
$2^j/4$ and length $4^j$ to each edge.  Then attach the removed disks
back on the far side of the cylinders.  This creates a simply
connected space with a trivial and thus bounded covering spectrum.

If we rescale $M$ by $1/2^j$ we get a tangent cone at infinity
which is not even semilocally simply connected.  It is a plane
with disks of radius $2^j/4$ centered at $(2^j,0)$ removed and
half cylinders attached for all values of $j\in \Bbb{Z}$.  Its
cut-off covering spectrum is clearly still trivial but its covering
spectrum is very large.

Note that without much difficulty, we could smooth $M$ to make it
a manifold and still get the same tangent cone at infinity.
\end{example}

\sect{Applications with Curvature Bounds} \label{sect6}

In this section we describe applications to complete
noncompact Riemannian manifolds
with lower bounds on their sectional and Ricci curvature
and their limit spaces.

\subsection{Sectional Curvature and the Soul Theorem}\label{sect6.1}

Cheeger-Gromoll \cite{Cheeger-Gromoll} proved that complete manifolds
with nonnegative
sectional curvature are diffeomorphic to normal bundles over totally
geodesic
compact submanifolds called souls.
Sharafutdinov \cite{Sharafutdinov} then
proved there
was a distance nonincreasing retraction to the soul: $P:M \to S$.
Perelman \cite{Perelman1, Perelman2} showed that
$P$ is a Riemannian submersion and extended the distance
nonincreasing retraction to complete Alexandrov  spaces
with nonnegative curvature.  Using the distance nonincreasing
retraction we can show that the covering spectrum of these
spaces behave exactly like the covering spectrum of a compact space.

\begin{theo} \label{sectthm}
If $M^n$ is a complete noncompact Alexandrov space with nonnegative
curvature, then
\be
CovSpec(M^n)=CovSpec(S^k)=CovSpec(\bar{T}_{R}(S^k))
\ee
 where $S^k$ is a soul
and $\bar{T}_R(S^k)$ is the $R$-closed tubular neighborhood around $S^k$.
\end{theo}

In light of the above paragraph Theorem~\ref{sectthm}
follows directly from the following theorem which we also use at the end of the 
paper.

\begin{theo} \label{alex}
If $M$ has a totally geodesic soul $S^k$
with a distance nonincreasing retraction $P: M\to S^k$
then
\be
CovSpec(M^n)=CovSpec(S^k)=CovSpec(\bar{T}_{R}(S^k)).
\ee
\end{theo}

Note that this is significantly stronger than the
loops to infinity property which says that curves
are homotopic outward.  In fact the curves are homotopic
inward to curves in the soul.

\Pf If $\delta \in CovSpec(M^n)$, then $\tilde{M}^{\delta'} \not=
\tilde{M}^\delta$ for all $\delta' > \delta$. Namely $\pi(M,\delta) \not =
\pi(M, \delta')$ for all $\delta' > \delta$. So for each $\delta' > \delta$,
there is $\gamma_{\delta'}$ in $\pi_1(M)$ such that $g_{\delta'}$ is
generated by elements lying $\delta'$-balls of $M$ but not generated by
elements lying $\delta$-balls of $M$. Since $P: M^n \ra S^k$ is distance
nonincreasing, $P$ maps balls of $M$ to the same or smaller size of ball of
$S^k$. Hence $P(\gamma_{\delta'})$ is generated by elements lying
$\delta'$-balls of $S^k$. Since $P$ is a retraction $P(\gamma_{\delta'})$ is
freely homotopic to $\gamma_{\delta'}$ so it can not be generated by
elements lying $\delta$-balls of $M$, therefore not $\delta$-balls of $S^k$.
Now for each $\delta' > \delta$, we have $P(\gamma_{\delta'})$ is generated
by elements lying $\delta'$-balls of $S^k$ but not $\delta$-balls of $S^k$. This
means $\delta \in CovSpec(S^k)$.

Conversely, if $\delta \in CovSpec(S^k)$, for each $\delta' > \delta$, there
is  $\gamma_{\delta'}$ in $\pi_1(S)$ such that $g_{\delta'}$ is generated by
elements lying $\delta'$-balls of $S$ but not generated by elements lying
$\delta$-balls of $S^k$. $g_{\delta'}$ is  not generated by elements lying
$\delta$-balls of $M^n$ either by above argument. Therefore $\delta \in
CovSpec(M^n)$.

Since $\pi_1(M^n) = \pi_1(S^k) = \pi_1(\bar{T}_{R}(S^k))$, we have
\be
CovSpec(\bar{T}_{R}(S^k)) \subset CovSpec(S^k)
\textrm{ and } CovSpec(M^n) \subset
CovSpec(\bar{T}_{R}(S^k)).
\ee
 Hence they are all equal.
\qed

\subsection{An Almost Soul Theorem}\label{sect6.2}

In this section we apply our results to complete
noncompact Riemannian manifolds with nonnegative sectional
curvature.  To do so we first study sequences of
manifolds with $sect \ge -\epsilon_i$ where $\epsilon_i$
converges to $0$ and prove an almost soul theorem:

\begin{theo} \label{almostsoul}
If $(X_i,x_i)$ are complete locally compact length spaces
converging in the pointed Gromov-Hausdorff sense
to a locally compact length space $(Y,y)$ such that
$Y$ is the normal bundle over a totally geodesic soul
with a distance nonincreasing retraction $P: Y \to S$
then there exist compact almost-souls $S_i \subset X_i$
with
\be
diam (S_i)=D_i \to diam(S) \textrm{ and }d_{X_i}(x_i,S_i) \to d_Y(y,S)
\ee
such that for any $b>a>0$ and any $R_2>R_1>0$ we have
\be
d_{H}(CovSpec(\bar{T}_{R_1}(S_i))\cap [a,b],
CovSpec(\bar{T}_{R_2}(S_i))\cap [a,b]) \to 0
\ee
where $\bar{T}_r(A)$ denotes the closed tubular neighborhood
about $A$ with the induced length metric.
\end{theo}

Note that the almost souls constructed here are not totally geodesic but
are compact. They are only soul like in the sense that loops
slide toward them, so that the covering spectrum is the same
on two distinct tubular neighborhoods.

\Pf
Note that $\bar{T}_{R_2}(S)$ is contained in some large
ball $B_{R_3}(y)$, and that there must be
an $\epsilon_i$ almost isometry $f_i: B_{x_i}(R_i) \to B_y(R_3)$.
Let $S_i$ be the closure of the
preimage of the soul $S\subset Y$:
\be
S_i=Cl(f_i^{-1}(S)).
\ee

Suppose the theorem is false. Then there exists
$b>a>0$ and $R_2>R_1>0$ and $\delta_i \in [a,b]$ such that
\be
\delta_i \in
CovSpec(\bar{T}_{R_1}(S_i))\setminus
CovSpec(\bar{T}_{R_2}(S_i))
\cup
CovSpec(\bar{T}_{R_1}(S_i))\setminus
CovSpec(\bar{T}_{R_2}(S_i)).
\ee
Since $\delta_i\in [a,b]$ a subsequence converges to
some $\delta \in [a,b]$.

By Theorem 8.4 of \cite{SoWei3}, applied to
$\bar{T}_{R_1}(S_i)$ and $\bar{T}_{R_2}(S_i)$
which converge to
$\bar{T}_{R_1}(S)$ and $\bar{T}_{R_2}(S)$, we
know
\be
\delta \in
CovSpec(\bar{T}_{R_1}(S))\setminus
CovSpec(\bar{T}_{R_2}(S))
\cup
CovSpec(\bar{T}_{R_1}(S))\setminus
CovSpec(\bar{T}_{R_2}(S)).
\ee
However no such $\delta$ exists by Theorem~\ref{alex}.
\qed

The following corollary follows immediately from
Theorem~\ref{almostsoul} and Theorem~\ref{sectthm}.

\begin{coro}\label{corosoul1}
Given any $h>0$, any $b>a>0$ and any $R_2>R_1>0$ there exists
$\epsilon=\epsilon(h,a,b,R_1,R_2)>0$ sufficiently small that
if $M^n$ has $sect \ge -\epsilon$ then
there is a compact $S \subset M^n$ such that
the Hausdorff distance:
\be
d_{H}(CovSpec(\bar{T}_{R_1}(S))\cap [a,b],
CovSpec(\bar{T}_{R_2}(S))\cap [a,b]) <h.
\ee
\end{coro}

Rescaling this corollary and consulting Theorem~\ref{almostsoul}
to locate the almost souls we get:

\begin{coro}\label{sectge-epsi}\label{corosoul2}
Given any $h,r,D>0$, any $b>a>0$ and any $R_2>R_1>0$ there exists
$\epsilon=\epsilon(h,a,b,R_1,R_2)>0$ sufficiently small that
if $M^n$ has $sect \ge - 1$ and $p\in M^n$
then there is a compact $S \subset M^n$ with $diam(S)\le D$
and $d(S,p)<r$
such that
\be \label{minisoul}
d_{H}(CovSpec(\bar{T}_{R_1\epsilon}(S))\cap [a\epsilon,b\epsilon],
CovSpec(\bar{T}_{R_2\epsilon}(S))\cap [a\epsilon,b\epsilon])
<h\epsilon.
\ee
\end{coro}

We can call such sets $S$ satisfying (\ref{minisoul})
subscaled souls and manifolds with this property
manifolds with many subscaled souls.

Note that a single space with thinner and thinner cylindrical subsets
would satisfy this corollary but a space with tiny handles would not.
A hyperbolic manifold will not have arbitrarily small handles,
but rather either looks locally thick like Euclidean space, or thin
like in a cusp where is it somewhat cylindrical.  Intuitively, this
corollary is saying manifolds with a uniform lower bound on sectional
curvature have a similar behavior.

\subsection{Nonnegative Ricci Curvature} \label{sect6.3}

When a complete noncompact manifold has nonnegative Ricci
curvature then it doesn't always have a soul.  However the
first author proved that such a manifold wither has
the loops to infinity property [Defn~\ref{defloopstoinfty}]
or it is the flat normal bundle over a compact totally geodesic
soul \cite{So-loops} [Theorem 11].  An example where the
latter occurs is the infinite Moebius strip.  Notice that the double cover
of the infinite Moebius strip is a flat cylinder.  In fact
Theorem 11 of \cite{So-loops} states that a double cover
always splits isometrically when the space fails to have the
loops to infinity property.  This means that is is the isometric
product of a line with another metric space.

This has profound implications on the cut-off covering spectrum:

\begin{theo}\label{Riccicutoff}
Let  $X$ be a complete noncompact manifold with $Ricci \ge 0$.
If $CovSpec_{cut}(X)$ is not empty then
$X$ has a double cover which splits isometrically
and $X$ is a flat normal bundle over a compact totally
geodesic soul, in which case $X$ has one element in
the covering spectrum and it is half the length of the shortest
closed geodesic which lifts as an open curve to this double cover.
\end{theo}

\begin{coro} \label{Riccicutcor}
If $M$ is a Riemannian manifold with Ricci curvature
strictly positive at one point
and $Ricci \ge 0$ everywhere then $CovSpec_{cut}(M)=\emptyset$.
\end{coro}

\noindent{\bf Proof of Theorem~\ref{Riccicutoff}:}
In \cite{So-loops} it is proven that a complete noncompact manifold, $M^m$,
with nonnegative Ricci curvature has the loops to infinity property
unless a double cover splits isometrically and $M^m$ is the flat normal
bundle over a compact totally geodesic soul {Theorem 7, Theorem 11].

When $M$ has
the loops to infinity property, we just apply
Theorem~\ref{thmloopstoinfty}.

When $M$ has a soul $S$ and
a split double cover, $\Bbb{R}\times K$, there is a collection
of loops $C$ which lift open to this double cover.
Each $C$ lifts to a curve of the form $(a,b)$ where
$a$ is a loop in $\Bbb{R}$ and $b$ is a loop in, $K$, the
compact double cover of the soul.  Note that any
curve $C$ is freely homotopic to the projection $\bar{C}$ of
the loop $(0,b)$ and $\bar{C}$ is shorter than $C$.

Let $\delta_0=\inf L(C)/2 =\inf L(\bar{C})/2$.  Then
$2\delta_0$ is the length of the shortest open path in $K$
which projects to a loop in the soul.  So it is positive
and is achieved by a closed geodesic $\gamma$ with represents some
element $g\in \pi_1(S,\gamma(0))\subset \pi_1(M,\gamma(0))$.

We claim $CovSpec_{cut}(X)=\{\delta_0\}$.

Given any closed curve $\sigma$ based at $\gamma(0)$
either $\sigma$ lifts to an open path or to a closed
loop in the double cover.  If $\sigma$ lifts as a closed
loop to the split double cover, then its lift has the
loops to infinity property, so we can project the homotopy down
and see that $\sigma$ has the loops to infinity property
and will not contribute to the cut-off covering spectrum.

So suppose $\sigma$ lifts as an open path to $\tilde{X}^{\delta,R}_{cut}$
for all values of $\delta$ and $R$.
Since the split double cover is a double cover,
$\sigma$ must be homotopic to $\gamma$ following
a loop which lifts to a closed loop to the double cover.
Since all loops in the split double cover have the loops to infinity
property,
that loop must as well.
So $\sigma$ must be $\delta, R$ homotopic to $\gamma$
for any value of $\delta$ and $R$.

Thus $\gamma$ alone
suffices to detect the distinct $\tilde{X}^{\delta,R}_{cut}$.
So there is only one element in the $R$ cut-off covering
spectrum for any value of $R$.  This element must be
$\delta_0$ because $\gamma$ was the shortest loop in
the class that doesn't  have the loops to
infinity property.
\qed

\subsection{Further Directions}  \label{sect6.4}

Here we discuss potential applications to manifolds with
$Ricci \ge -(n-1)$.  These applications will arise by studying
limits of manifolds with $Ricci \ge -\epsilon_i$.  These
limits spaces are complete locally compact length spaces \cite{Gr}.

Let $X$ be the Gromov-Hausdorff limit of complete noncompact
Riemannian manifolds $M_i$ with $Ricci(M_i)\ge -\epsilon_i \to 0$.
Such a length space is very similar to a Riemannian manifold with
nonnegative Ricci curvature.  Cheeger and Colding have proven
that the splitting theorem holds on such a space \cite{ChCo1}
and the authors have proven that such spaces have universal covers
\cite{SoWei2}.  Combining these facts with the proof of
Theorem 7 in \cite{So-loops}, the authors proved
that $X$ either has the loops to infinity property or
the universal cover splits isometrically \cite{SoWei2}[Cor 4.9].

Thus, using the proof of Theorem~\ref{Riccicutoff}, we can conclude
the following:

\begin{theo}\label{Riccicutoff2}
If $X$ is a limit space as described above, then either
$CovSpec_{cut}(X)$ is empty or its universal cover
splits isometrically.
\end{theo}

Now Theorem 11 in \cite{So-loops}, which states that when $X$
doesn't have the loops to infinity property then $X$ has a split
double cover and $X$ has a compact soul, was never extended.
Its proof involves differentiation.  In light of Cheeger's recent
work on differentiability of metric measure spaces, we conjecture
that the full theorem holds:

\begin{conj} \label{Conj1}
If $X$ is a limit space as in the above paragraph, then either
$X$ has the loops to infinity property or it has a split double cover and
it is the flat normal bundle over a compact totally geodesic soul.
\end{conj}

If Conjecture~\ref{Conj1} holds then the proof of Theorem~\ref{Riccicutoff}
extends to that setting and the following conjecture holds:

\begin{conj} \label{conj2}
Theorem~\ref{Riccicutoff} holds for such limit spaces $X$.
\end{conj}

Conjecture~\ref{conj2} is of particular interest because it has implications
to manifolds with $Ricci \ge -(n-1)$.  In particular, the application
looks something like a Margulis lemma but with distinct implications:

\begin{conj} \label{conj3}
Given a complete Riemannian
$M^n$ with $Ricci \ge -1$ or a Gromov-Hausdorff limit
of such spaces, for all $b>a>0$ there exists
$\rho=\rho(a,b,n)$ such that
\be
CovSpec_{cut}^{\rho}(M,p)\cap [a\rho,b\rho]=\emptyset
\ee
or there is a subscaled soul as in (\ref{minisoul})
near $p$.
\end{conj}

Note that $\rho$ does not depend on the manifold or the basepoint.
It is essentially saying that the manifold looks locally like
a manifold with $sect \ge -1$ in neighborhoods where loops don't
slide outward.  One might think of this as saying that most small loops
on the manifold slide around, and when you hit a location where
they don't slide there is a kind of twisting effect similar
to a Moebius strip.

Returning from the intuitive to the concrete, we can prove:

\begin{theo}\label{conj2to3}
Conjecture~\ref{conj2} implies Conjecture~\ref{conj3}
\end{theo}

\noindent{\bf Proof of Theorem~\ref{conj2to3}}
Suppose $(M_i,d_i)$ have $Ricci \ge -1$ and $\rho_i \to 0$
with
\be
\delta_i \rho_i \in CovSpec_{cut}^{\rho_i}(M_i,d_i,p_i),
\ee
where $\delta_i \in [a,b]$.
Rescaling $M_i$ by $rho_i$, Lemma~\ref{rescale} says
\be
\delta_i  \in CovSpec_{cut}^{1}(M_i, d_i/\rho_i,p_i).
\ee
Gromov's compactness theorem implies that a subsequence
of the $(M_i, p_i, d_i/\rho_i)$ converge to some $(Y,y,d)$ which satisfies
the conditions of Theorem~\ref{Riccicutoff}. Thus
$CovSpec_{cut}(Y)=\emptyset$ or a split double cover.

Taking a further subsequence we can guarantee $\delta_i \to \delta\in
[a,b]$.
By Theorem~\ref{covcutconv} (\ref{ccca}) 
\be
\delta  \in CovSpec_{cut}^{1}(Y, d, y) \subset CovSpec_{cut}(Y).
\ee
Thus $Y$ has a split double cover and is the flat normal
bundle over a compact totally geodesic soul $S$.

Finally we apply Theorem~\ref{almostsoul} to $Y$.
\qed


\sect{Appendix A} \label{AxA}

As the concept of lower semiclosure does not seem to
appear in the literature, we include a brief exposition
here.

\begin{defn} \label{defnlowersemiclosed}
 A lower semiclosed subset of the
real line is a set $A$ such that $lim_{j\to\infty} a_j \in A$
whenever
$a_j$ is a decresing sequence of elements of $A$.
\end{defn}

\begin{defn} \label{defnlowersemiclosure}
The lower semiclosure of a set $A$, denoted
$Cl_{lower}(A)$, is the intersection
of all lower semiclosed sets containing $A$.
\end{defn}

\begin{lemma}
The lower semiclosure of $A$ is the union of $A$ and the
limits of any decreasing sequence of $a_j\in A$.
\end{lemma}

\begin{lemma}\label{appendixopen}
If $x\notin A$ and $A$ is lower semiclosed, then there
exists $\epsilon>0$ such that
\be
[x,x+\epsilon)\cap A =\emptyset.
\ee
\end{lemma}

The following theorem implies that
$CovSpec(X)$ and $CovSpec_{cut}(X)$ 
are lower semiclosed subsets of $\Bbb(R)$.

\begin{theo} \label{appendixthm}
Let $X_s$ be a collection of metric spaces parametrized by
a real line, $s\in \Bbb(R)$, such that whenever
$s_1<s_2$ we have $X_{s_1}$ covers $X_{s_2}$.
Any set $A$ defined as follows:
\be
A := \{ s: \forall s' > s X_{s'}\neq X_{s}\}
\ee
then $A$ is lower semiclosed.
\end{theo}

\Pf
Let $s_j \in A$ a decreasing sequence converging to $s_\infty$.
We need to show $s_\infty \in A$.  Let $s'>s_\infty$.  Then
for $j$ sufficiently large, we have $s'>s_j$.  Since
$s_j\in A$ this means $X_{s'} \neq X_{s_j}$ so
$X_{s_j}$ is a nontrivial cover of $X_{s'}$.
And since $X_{s_\infty}$ covers $X_{s_j}$, it
must be a nontrivial cover of $X_{s'}$ as well.
\qed

\begin{example}
If $A_j$ are all lower semiclosed sets, the $\bigcup_{j\in \Bbb{N}}A_j$
need not be lower semiclosed.  For example, let
$A_j=\{1+1/k: k=1,2,...j\}$.
\end{example}

\sect{Appendix B} \label{AxB}

This appendix provides a minor correction to \cite{SoWei3}
adding a hypothesis to Lemma 5.8 [Lemma~\ref{balltoL1}] and
proving a new related Lemma~\ref{balltoL} which circumvents
the additional hypothesis in Theorem 5.7 of  \cite{SoWei3}
which applied Lemma 5.8.  
Both Lemmas~\ref{balltoL1} and~\ref{balltoL} are applied in this
paper as well.

In \cite{SoWei3}, there is an omission in the statement of
Lemma 5.8.  The proof requires that the curve
be rectifiable as pointed out to us by Conrad Plaut.
The corrected statement is:

\begin{lem} \label{balltoL1} 
Given a complete length space $Y$, and suppose
 $C:[0,L] \to B_q(\delta)\subset Y$ is rectifiable 
 then $C$ is freely homotopic in $B_q(\delta)$
 to a product of curves, $C_i$, of length $L(C_i)<2\delta$ 
all of which lie in $B_q(\delta)$.
 \end{lem}

The proof is as in \cite{SoWei3}:

\noindent{\bf Proof of Lemma~\ref{balltoL1}:}
Since $B_q(\delta)$ is open and the image of $C$ is closed
there exists $\epsilon>0$ such that $Im(C) \subset B_q(\delta-\epsilon)$.
Take a partition, $0=t_0<t_1<...<t_k=L$, such that $t_{j+1}-t_j<\epsilon$,
and let $\gamma_j$ run minimally from $q$ to $C(t_j)$ making
sure to choose $\gamma_0=\gamma_k$.    Then
$C$ is clearly freely homotopic in $B_q(\delta)$ to the combination
$\gamma_j C([t_j,t_{j+1}]) \gamma_j^{-1}$, and each of these
curves has length $< 2(\delta-\epsilon)+\epsilon<2\delta$.
\qed

This lemma was applied to prove Theorem 5.7 of \cite{SoWei3}
that the covering spectrum is determined by the marked length
spectrum when $X$ is a compact length space.   In this Appendix we
provide a correction of that proof which clarifies why
we can select a rectifiable curve before applying
the corrected Lemma 5.8.  

It should be noted that the Hawaii Ring with loops of length $1/n$
has a nonrectifiable curve $C$ which traverses all of its loops
by traveling faster and faster.  So not every homotopy class of
curves in a complete length space contains 
a rectifiable representative.  

Recall that in the definition of the delta covers and
$\pi(Y,\delta,p)$, curves $\beta$ are classified according to
their location: the fact that they are contained in a ball
of radius $\delta$.  To prove Theorem 5.7, we showed that
we could control the lengths of representatives curves as well
using Lemma~\ref{balltoL1}.  This relationship is correctly
stated here and we provide a proof which clarifies how we
can select a rectifiable curve:

\begin{lemma} \label{balltoL} 
Given a complete length space $Y$, and suppose
$\delta'<\delta$ and $\tilde{Y}^\delta \neq \tilde{Y}^{\delta'}$
 or equivalently, $\pi(Y,\delta,p)\neq \pi(Y, \delta',p)$,
then there is a 
rectifiable curve $\beta$ of length
$L(\beta)<2\delta$ and a curve $\alpha$ running from $p$
to $\beta(0)$, such that
\be \label{lookherea}
[\alpha^{-1}\circ\beta\circ\alpha] \notin \pi(Y,\delta',p).
\ee
So $\beta$ lifts closed to $\tilde{X}^\delta$ and
open to $\tilde{X}^{\delta'}$.
\end{lemma}

Note that by the definitions, we already know there exists
a $\beta '$ satisfying (\ref{lookherea}) whose image lies in 
a ball of radius $\delta$, otherwise all the generators of $\pi(Y,\delta,p)$
would already lie in $\pi(Y, \delta',p)$ and $\tilde{Y}^\delta$
would equal $\tilde{Y}^{\delta'}$.  The
difficult part is proving that we can control its length.
If we knew such a $\beta'$ were rectifiable, then we
could apply Lemma~\ref{balltoL1} with $\beta '=C$ to get: 
$$
[\alpha^{-1}\circ\beta '\circ\alpha]= [\alpha^{-1}\circ C_1\circ\alpha]\cdot
 [\alpha^{-1}\circ C_2\circ\alpha]\cdots [\alpha^{-1}\circ C_N\circ\alpha].
$$
We could then select one of the $C_i$ to be $\beta$
and (\ref{lookherea}) would be satisfied
since otherwise our original curve would be in $\pi(Y,\delta',p)$.  

When $\beta'$ is not rectifiable, as can occur on the Hawaii Ring,
then we need only replace it with a rectifiable curve that has the same
lifting properties.  

The following proof of Lemma~\ref{balltoL} begins with a technique
suggested by Conrad Plaut to shift a nonrectifiable curve
to one with similar lifting properties.

\vspace{.2cm}

\noindent {\bf Proof of Lemma~\ref{balltoL}}:
First we know there is a continuous curve
$\beta '$ contained in some $B_q(\delta)$
which satisfies (\ref{lookherea}).  
Note that since $[0,1]$ is
compact, the image of $\beta '$ is a compact set.  So there
exists an $\epsilon\in (0, \delta'/10)$ sufficiently small that 
for any $t\in S^1$, $B_{\beta '(t)}(5\epsilon)\subset B_q(\delta)$.
By continuity, we can create a partition $0=t_0<t_1<t_2<...<t_N=1$
so that each segment $\beta([t_i,t_{i+1}])$ lies in one of
these balls.

Note that a curve, $\eta_i$, created by running from $\beta '(t_i)$ to
$\beta '(t_{i+1})$ along $\beta '([t_i, t_{i+1}])$ and then
back to $\beta '(t_i)$ along a minimal geodesic, lies
in a ball of radius $4\epsilon < \delta'$.  Choosing appropriate
$\alpha_i$ running along $\alpha$ and then up $\beta '([0,t_i])$
we see that 
\be
[\alpha_i^{-1}\circ\eta_i\circ\alpha_i] \in \pi(Y,\delta',p).
\ee
So we can construct a rectifiable curve
$\beta$ using piecewise minimizing geodesics
running between $\beta '(t_i)$.  This $\beta$ will lie
in $B_q(\delta)$, it will be rectifiable, and it will satisfy
(\ref{lookherea}) since $\beta'$ did but the $\eta_i$ do not.  

We now apply Lemma~\ref{balltoL1} to this curve $\beta$, and
we see it is generated by $\beta_i$ of length $\le 2\delta$
contained in the same $B_q(\delta)$.  At least one of these
curves must satisfy (\ref{lookherea}), or $\beta$ would not.
\qed

Theorem 5.7 of \cite{SoWei3} states that on a compact length
space with a universal cover
the marked length spectrum determines the covering spectrum.
This proof is now corrected applying Lemma~\ref{balltoL}
in place of Lemma 5.8 to show that we can use rectifiable
representative loops of length $<2\delta'$ for the $\beta_i$
near the end of the proof of the theorem.  Note that 
Theorem~\ref{balltoloop} is essentially contains a new proof 
of the same theorem in simpler language.

 \end{document}